\documentclass{article}

\usepackage{lineno}
\usepackage{hyperref}

\usepackage[top=4cm, bottom=4cm, left=3cm, right=3cm]{geometry}
\usepackage{amsmath} 
\usepackage{amsthm} 
\usepackage{mathrsfs}
\usepackage{amsfonts}
\usepackage{mathtools}
\usepackage{amssymb}
\usepackage{xcolor}
\usepackage{tikz-cd}

\tikzset{
	symbol/.style={
		draw=none,
		every to/.append style={
			edge node={node [sloped, allow upside down, auto=false]{$#1$}}}
	}
}

\hypersetup{pdfpagemode={UseOutlines},	bookmarksopen, pdfstartview={FitH},	colorlinks,	linkcolor={red},	citecolor={blue},	urlcolor={blue}}

\newtheorem{thm}{Theorem}[section]

\newtheorem{proposition}[thm]{Proposition}
\newtheorem{theorem}[thm]{Theorem}
\newtheorem{lemma}[thm]{Lemma}

\theoremstyle{definition}
\newtheorem{definition}[thm]{Definition}

\theoremstyle{remark}
\newtheorem{remark}[thm]{Remark}
\numberwithin{equation}{section}

\newcommand{\RR}{\mathbb{{R}}}

\newcommand{\NN}{\mathbb{{N}}}

\newcommand{\CC}{\mathbb{{C}}}
\newcommand{\DD}{\mathbb{{D}}}
\newcommand{\nul}{\mbox{\normalfont{nul}}}
\newcommand{\defec}{\mbox{\normalfont{def}}}
\newcommand{\ascent}{\alpha}
\newcommand{\descent}{\delta}

\newcommand{\dom}{\mathcal{D}}
\newcommand{\kernel}{\mathcal{N}}
\newcommand{\ran}{\mathcal{R}}

\newcommand{\codim}{\mbox{\normalfont{codim}}}
\newcommand{\linearOp}{\mathcal{L}}
\newcommand{\compact}{\mathcal{K}}

\newcommand{\BSect}{\text{\normalfont{BSect}}}

\title{Spectral mapping theorems for essential spectra and regularized functional calculi}
\author{Jesús Oliva-Maza}

\begin{document}

\maketitle

\begin{abstract}
Gramsch and Lay \cite{gramsch1971spectral} gave spectral mapping theorems for the Dunford-Taylor calculus of a closed linear operator $T$,
$$\widetilde{\sigma}_i(f(T)) = f(\widetilde{\sigma}_i(T)),
$$
for several extended essential spectra $\widetilde{\sigma}_i$. In this work, we extend such theorems for the regularized functional calculus introduced by Haase \cite{haase2005general, haase2005spectral} assuming suitable conditions on $f$. At the same time, we answer in the positive a question posed by Haase \cite[Remark 5.4]{haase2005spectral} regarding the conditions on $f$ which are sufficient to obtain the spectral mapping theorem for the usual extended spectrum $\widetilde \sigma$. We use the model case of bisectorial-like operators, although the proofs presented here are generic, and are valid for similar functional calculi.
\end{abstract}




\section{Introduction}

Let $T$ be a bounded self-adjoint operator on a Hilbert space $\mathcal H$. A spectral singularity $\lambda \in \sigma(T)$, where $\sigma(T)$ denotes the spectrum of $T$, is said to be in the essential spectrum of $T$ is $\lambda$ is not an isolated eigenvalue of finite multiplicity, see \cite{wolf1959essential}.

If $T$ is not self-adjoint, or if $T$ is an operator on a Banach space $X$, most modern texts define the essential spectrum $\sigma_{ess}(T)$ of $T$ in terms of Fredholm operators, that is, $\lambda \in \sigma_{ess}(T)$ iff $\lambda - T$ is not a Fredholm operator. Recall that a bounded operator $T$ is a Fredholm operator if both its nullity $\nul(T)$ (dimension of its kernel) and its defect $\defec(T)$ (codimension of its range) are finite. One of the main useful properties of the essential spectrum (defined this way) is that it is invariant under compact perturbations. As a matter of fact, $\sigma_{ess}(T) = \sigma(p(T))$, where $p(T)$ is the projection of $T$ in the Calkin algebra, i.e., the quotient algebra of the bounded operators $\linearOp(X)$ on $X$ modulo the compact operators $\compact(X)$ on $X$.

However, several different definitions for the essential spectrum were introduced in the 50s and 60s, especially in the framework of differential operators. For instance, if we denote by $G^l$ ($G^r$) the semigroup of left (right) regular elements in the Calkin algebra $\linearOp(X)/\compact(X)$, Yood \cite{yood1951properties} studied the spectral sets of $T\in \linearOp(X)$ given by $\sigma_2(T):=\{\lambda \in \CC \, : \, \lambda - p(T) \notin G^l\}$ and $\sigma_3(T) := \{\lambda \in \CC \, : \, \lambda - p(T) \notin G^r\}$. Indeed, he obtained the following characterizations of the semigroups $p^{-1}(G^l), p^{-1}(G^r)$:
\begin{align*}
	p^{-1} (G^l) &= \{T \in \linearOp(X) \, : \, \nul(T) < \infty \mbox{ and } \ran(T) \mbox{ is complemented}\}, \\
	p^{-1} (G^r) &= \{T \in \linearOp(X) \, : \, \defec(T) < \infty \mbox{ and } \kernel(T) \mbox{ is complemented}\},
\end{align*}
where $\ran(T), \, \kernel(T)$ denote the range space of $T$ and the kernel of $T$ respectively.

Alternatively, spectral sets associated with semi-Fredholm operators have also been referred as essential spectra. More precisely, let $\Phi^+$, $\Phi^-$ be given by
\begin{align*}
	\Phi^- &:= \{T \in \linearOp(X) \, : \, \nul(T) < \infty \mbox{ and } \ran(T) \mbox{ is closed}\},
	\\\Phi^+ &:= \{T \in \linearOp(X) \, : \, \defec(T) < \infty\}.
\end{align*}
Then Gustafson and Weidmann \cite{gustafson1969essential} used the term essential spectra for the spectral sets $\sigma_4(T) := \{\lambda \in \CC \, : \, \lambda - T \notin \Phi^-\}$, $\sigma_5(T) := \{\lambda \in \CC \, : \, \lambda - T \notin \Phi^+\}$, while Kato \cite{kato1966perturbation} considered the spectral set $\sigma_6(T):= \sigma_4(T) \cap \sigma_5(T)$, i.e., $\lambda \in \sigma_6(T)$ iff $\lambda - T$ is not in $\Phi^- \cup \Phi^+$. Note that $\sigma_2(T) = \sigma_4(T)$ and $\sigma_3(T) = \sigma_5(T)$ if $T$ is an operator on a Hilbert space, but these equalities are not true in general in the framework of Banach spaces, see the work of Pietsch \cite{pietsch1960theorie}.

In another direction, Browder \cite{browder1961spectral} defined the essential spectrum of $T$ ($\sigma_8(T)$ here) as those spectral values of $T$ which are not isolated eigenvalues of finite multiplicity of $T$ nor isolated eigenvalues of finite multiplicity of the adjoint operator $T^\ast$ (cf. \cite{lay1968characterizations}). It turns out that $\lambda \notin \sigma_8(T)$ iff $\lambda$ is a pole of the resolvent of finite rank \cite[Lemma 17]{browder1961spectral}. With this in mind, Gramsch and Lay \cite{gramsch1971spectral} considered the following another essential spectrum, which is given by
$$\sigma_9(T) := \{\lambda \in \CC \, : \, \mbox{the resolvent of } T \mbox{ is not meromorphic at } \lambda\}.
$$
Nevertheless, the essential spectrum of Browder $\sigma_8(T)$ fails to satisfy the property of invariance under compact perturbations. In this regard, Schechter \cite{schechter1966essential} defined the essential spectrum of $T$, $\sigma_7(T)$ here, as the larger subset of $\sigma(T)$ which is invariant under compact perturbations. Equivalently $\lambda \notin \sigma_7(T)$ iff $\lambda - T$ is Fredholm with index zero, i.e., $\nul(\lambda - T) = \defec(\lambda - T) < \infty$.

In this work, we deal with spectral mapping theorems for the different essential spectra described above, that is, identities of the form
\begin{equation}\label{spectralMappingTheoremIntro}
	\sigma_i(f(T)) = f (\sigma_i(T)).
\end{equation}
There, $f$ is a function in the domain of a functional calculus of a (possibly unbounded) operator $T$.

On this point, the first approach to a Banach space functional calculus of unbounded operators is the so-called Dunford-Taylor calculus. For this calculus, one considers functions $f$ which are holomorphic in an open set containing the extended spectrum $\widetilde \sigma(T)$ of $T$, defined by $\widetilde \sigma(T) := \sigma(T) \cup \{\infty\}$ if $T$ is unbounded and $\widetilde \sigma(T) := \sigma(T)$ otherwise. Then, $f(T)$ is determined by
\begin{equation}\label{dunfordTaylorFormula}
	f(T) := f(\infty) + \frac{1}{2\pi i} \int_\Gamma f(z) (z-T)^{-1} \, dz,
\end{equation}
where $\Gamma$ is a suitable finite cycle that avoids $\widetilde\sigma(T)$. Moreover, the Dunford-Taylor formula above \eqref{dunfordTaylorFormula} still works if the curve $\Gamma$ touches $\widetilde\sigma(T)$ at some points $a_1,\ldots,a_n$ and $f$ is not holomorphic at $a_1, \ldots,a_n$, as long as $f$ tends to a finite number at each point $a_1,\ldots,a_n$ fast enough to deal with the size of the resolvent at these points. In this case, we say that $f$ has regular limits (at $a_1,\ldots,a_n$) and denote it by $\mathcal E(T)$.

Furthermore, in the setting of strip-type operators, Bade \cite{bade1953operational} introduced a `regularization trick' in order to define $f(T)$ for functions which do not grow too fast at $\infty$. This `regularization trick' was further developed, in the framework of sectorial operators, by McIntosh \cite{mcintosh1986operators}, Cowling et al \cite{cowling1996banach} and Haase \cite{haase2005general}. In particular, fractional powers and/or logarithms can be defined for suitable operators with this `regularization trick'.

Here, we consider the `regularized' functional calculus of meromorphic functions developed by Haase \cite{haase2005general}, which is based on the following idea. A meromorphic function $f$ is the domain of the regularized functional calculus of $T$, which we denote by $f \in \mathcal M(T)$, if there exists a holomorphic function $e \in \mathcal E(T)$ such that $e(T)$ is injective and $ef \in \mathcal E(T)$. In this case, one defines
\begin{equation}\label{haaseFunctCalc}
	f(T) := e(T)^{-1} (ef)(T),
\end{equation}
which is a (possibly unbounded) closed operator on $X$.

Going back to the spectral mapping theorem, identities \eqref{spectralMappingTheoremIntro} were proven by Gramsch and Lay \cite{gramsch1971spectral} in the setting of the Dunford-Taylor calculus, for most (extended) essential spectra described here, see Section \ref{functionalCalculusSection} for their definitions. Moreover, Gonz{\'a}lez and Onieva \cite{gonzalez1985spectral} used a unified approach and gave simpler proofs for these spectral mapping theorems. Their proofs are based on the following observations:
\begin{enumerate}
	\item[1)] a closed operator $T$ with non-empty resolvent set is essentially invertible iff the bounded operator $T(b-T)^{-1}$ is essentially invertible \cite[Lemma 1]{gonzalez1985spectral},
	\item[2)] for $f,g$ in the domain of the Dunford-Taylor calculus of $T$, one has $(fg)(T) =f(T)g(T)=g(T)f(T)$. As a consequence, $(fg)(T)$ is essentially invertible if/only if (see \cite[Lemma 3]{gonzalez1985spectral}) both $f(T),g(T)$ are essentially invertible,
\end{enumerate}
where we say that an operator $A$ is essentially invertible (regarding the essential spectrum $\sigma_i$) if $0 \notin \sigma_i(A)$, and
\begin{enumerate}
	\item[3)] if $f$ is in the domain of the Dunford-Taylor calculus of $T$, then one can assume that $f$ has a finite number of zeroes of finite multiplicity. 
\end{enumerate} 

It sounds sensible to ask whether these spectral mapping theorems can be extended to cover the functions in the domain of the regularized functional calculus given by \eqref{haaseFunctCalc}. This is partly motivated by potential applications in Fredholm theory, in particular if considering fractional powers or logarithms of unbounded operators. For instance, we make use, in an ongoing work with L. Abadías, of the results presented here to provide the essential spectrum of fractional Ces\`aro operators and H\"older operators acting on spaces of holomorphic functions. However, there are two main difficulties for such an extension of the spectral mapping theorem. First, for $f,g \in \mathcal M(T)$, it is not true in general that $(fg)(T) = f(T)g(T) =g(T)f(T)$, so item 2) above fails. Indeed, one only has the inclusions $f(T)g(T), g(T)f(T) \subseteq (fg)(T)$, where $S\subseteq T$ means that $\dom(S) \subseteq \dom(T)$ with $Sx = Tx$ for every $x \in \dom(S)$. Secondly, since the function $f$ may not be holomorphic at the points $a_1,\ldots,a_n$ where the integration path touches $\sigma(T)$, item 3) above also fails to be true.

Nevertheless, in the setting of sectorial operators, Haase \cite{haase2005spectral} overcame these two problems for the usual extended spectrum $\widetilde \sigma$, and obtained the spectral mapping theorem
\begin{equation}\label{haaseSpectralMT}
	\widetilde \sigma(f(T)) = f(\widetilde\sigma(T)),
\end{equation}
for a meromorphic function $f$ in the domain of the regularized functional calculus, i.e., $f \in \mathcal M(T)$, such that $f$ has almost logarithmic limits at the points $a_1,\ldots,a_n$ where the integration path $\Gamma$ touches $\widetilde\sigma(T)$. This `almost logarithmic' condition on the behavior of the limits of $f$ is stronger than asking $f$ to have regular limits at $a_1,\ldots,a_n$. As a matter of fact, Haase leaves open the question whether the hypothesis of $f$ having regular limits is sufficient to obtain the spectral mapping theorem, see \cite[Remark 5.4]{haase2005spectral}.

Still, it is far from trivial to extend the spectral mapping theorem \eqref{haaseSpectralMT} from the usual extended spectrum to the (extended) essential spectra described here. This extension, which is given here, is the main contribution of the paper. Even more, we obtain spectral mapping theorems for the essential spectra described above and for functions $f$ with regular limits lying in the domain of the regularized functional calculus of meromorphic functions \eqref{haaseFunctCalc}, answering in the positive Haase's conjecture on regular limits explained above.

To obtain these results, on the one hand, we provide a slightly simpler proof for the spectral inclusion of the usual extended spectrum, i.e., $f(\widetilde \sigma(A)) \subseteq \widetilde \sigma(f(A))$, than the one given in \cite{haase2005spectral}. As a matter of fact, we no longer make use of the composition rule of the functional calculus, which in the end allow us to weaken the condition on the function $f$ from almost logarithmic limits to the (quasi-)regular limits, cf. \cite[Remark 5.4]{haase2005spectral}. 

On the other hand, we address the items 2) and 3) above to cover all the essential spectra described above. First, we provide a commutativity property in Lemma \ref{domainLemma}, which is a refinement of \cite[Lemma 4.2]{haase2005spectral}, and which is crucial to deal with regularized functions $f \in \mathcal M(T)$ which are not in $\mathcal E(T)$. Moreover, since item 3) is not true if $f$ has a zero at the points $a_1,\ldots,a_n$ of $\widetilde \sigma(T)$ (recall that $f$ may not be meromorphic there), another issue of importance is whether the points $f(a_1),\ldots,f(a_n)$ belong to $\widetilde \sigma_i(f(T))$. To solve this, we apply different techniques depending on the topological properties (relative to $\widetilde\sigma(T)$) of $a_1,\ldots,a_n$. If these points are isolated points of $\widetilde \sigma(T)$, we provide useful properties of the spectral projections associated with such points in Lemmas \ref{resolventLemma} and \ref{quotientNFCLemma}. If otherwise, these points are limit points of $\widetilde \sigma(T)$, we make use of a mixture of topological properties shared by all the essential spectra considered here, and a mixture of algebraic properties of the regularized functional calculus in Propositions \ref{essentialInclusionProp1} and \ref{essentialInclusionProp2}.

In this work, we use the model case of bisectorial-like operators, which is a family of operators that slightly generalizes the one of bisectorial operators, see for instance \cite{arendt2006maximal,mielke1987maximalel}. This is partly motivated by two reasons. On the one hand, we want our results to cover the case when $T$ is the generator of an exponentially bounded group. This is because, in a forthcoming paper, we obtain spectral properties of certain integral operators via subordination of such operators in terms of an exponentially bounded group, namely, a weighted composition group of hyperbolic symbol. On the other hand, the another incentive to do this is the fact that the regularized functional calculus for bisectorial-like operators is easily constructed by mimicking the regularized functional calculus of sectorial operators \cite{haase2005general, haase2006functional}. Finally, bisectorial operators play an important role in the field of abstract inhomogeneous differential equations over the real line, so we are confident that our results have applications of interest in that topic.

Nevertheless, the proofs presented here are generic and are valid for every regularized functional calculus of meromorphic functions (in the sense of \cite{haase2005general, haase2006functional}) satisfying the properties collected in Lemmas \ref{functionalCalculusProperties}, \ref{pointSpectraLemma}, \ref{regularizerLemma} and \ref{differentNFCLemma}. For instance, our proofs work for the regularized functional calculus of sectorial operators and the regularized functional calculus of strip-type operators, as we point out in Subsection \ref{sectorialSubSect}.

The paper is organized as follows. The regularized functional calculus for bisectorial-like operators is detailed in Section \ref{functionalCalculusSection}. In Section \ref{noSingularPointsSection}, we give the spectral mapping theorems for a bisectorial-like operator $A$ with no singular points in the case the integration path $\Gamma$ does not touch any point of $\widetilde \sigma(T)$. The general case is dealt with in Section \ref{spectralMappingSection}. We give some final remarks in Section \ref{finalRemarksSection}, such as the answer in the positive to Haase's conjecture \cite[Remark 5.4]{haase2005spectral}.

\section{Extended essential spectra and regularized functional calculus for bisectorial-like operators}\label{functionalCalculusSection}

Let us fix (and recall) the notation through the paper. $X$ will denote an infinite dimensional complex Banach space. Let $\mathcal{L}(X)$, $C(X)$ denote the sets of bounded operators and closed operators on $X$, respectively. For $T \in C(X)$, let $\dom(T), \, \ran(T), \, \kernel(T)$ denote the domain, range, null space of $T$, respectively. Moreover, we denote the nullity of $T$ by $\nul(T)$, and the defect of $T$ by $\defec(T)$. The \textit{ascent} of $T$, $\ascent(T)$, is the smallest integer $n$ such that $\kernel(T^n) = \kernel(T^{n+1})$, and the \textit{descent} of $T$, $\descent(T)$, is the smallest integer $n$ such that $\ran(T^n) = \ran(T^{n+1})$. 

Now we recall the definition of the different essential spectra described in the Introduction. Following the notation and terminology of \cite{gonzalez1985spectral,gramsch1971spectral}, set
\begin{align*}
	\Phi_0 &:= \{T \in C(X) \, | \, \nul(T) = \defec(T) = 0\},
	\\ \Phi_1 &:= \{T \in C(X) \, | \, \nul(T), \defec(T) < \infty\},
	\\ \Phi_2 &:= \{T \in C(X) \, | \, \nul(T) < \infty, \, \ran(T) \mbox{ complemented}\},
	\\ \Phi_3 &:= \{T \in C(X) \, | \, \defec(T) < \infty, \, \kernel(T) \mbox{ complemented}\},
	\\ \Phi_4 &:= \{T \in C(X) \, | \, \nul(T) < \infty, \, \ran(T) \mbox{ closed}\},
	\\ \Phi_5 &:= \{T \in C(X) \, | \, \defec(T) < \infty\},
	\\ \Phi_6 &:= \Phi_4 \cup \Phi_5,
	\\ \Phi_7 &:= \{T \in C(X) \, | \, \nul(T) = \defec(T) < \infty\},
	\\ \Phi_8 &:= \{T \in \Phi_7 \, | \, \ascent(T) = \descent(T) < \infty\},
	\\ \Phi_9 &:= \{T \in C(X) \, | \, \ascent(T), \descent(T) < \infty\}.
\end{align*}
We observe that these operator families satisfy the following spectral inclusions

\begin{equation*}
\begin{tikzcd}[row sep = 0em, column sep=.7em]
	&&&& \Phi_3 \arrow[r,symbol=\subseteq] & \Phi_5 \arrow[rd, symbol=\subseteq] \\
	\Phi_0 \arrow[r,symbol=\subseteq] & \Phi_8 \arrow[r,symbol=\subseteq] &\Phi_7 \arrow[r,symbol=\subseteq] &\Phi_1 \arrow[ru,symbol=\subseteq] \arrow[rd, symbol=\subseteq] &&& \Phi_6 \qquad \mbox{and}
	\quad & \Phi_0 \arrow[r,symbol=\subseteq] & \Phi_8 \arrow[r,symbol=\subseteq]&\Phi_9,
	 \\
	&&&& \Phi_2 \arrow[r,symbol=\subseteq] & \Phi_4 \arrow[ru, symbol=\subseteq]
\end{tikzcd}
\end{equation*}
Then, the respective spectra $\sigma_i(T)$ are defined in terms of the above families by
$$\sigma_i(T) := \{\lambda \in \CC \, | \, \lambda - T \notin \Phi_i\} \quad \mbox{for} \quad i \in  \{0,1,...,9\}.
$$ 
Note that $\sigma_0(T)$ is the usual spectrum $\sigma(T)$ and most modern text use the term essential spectrum to denote the set $\sigma_1(T)$. It is also worth saying that the works \cite{gonzalez1985spectral, gramsch1971spectral}  also considered the essential spectrum $\sigma_{10}(T)$ defined in terms of normally solve operators, i.e., operators with closed range, see \cite{dunford1963linear}. However, there exist bounded operators on Hilbert spaces for which $\sigma_{10}(T^2) \not\subseteq (\sigma_{10}(T))^2$ and $\sigma_{10}(S^2) \not\supseteq (\sigma_{10}(S))^2$, see \cite[Section 5]{gramsch1971spectral}.

Next we define the extended essential spectra $\widetilde{\sigma}_i(T)$.
\begin{definition}\label{extendedEssSpectraDef}
	Let $T \in C(X)$. We define
	\begin{align*}
		\widetilde{\sigma}_i (T) &:= \begin{cases}
			\sigma_i(T) \mbox{ if } \begin{cases}
				\dom(T) = X, \, \mbox{ for } i \in \{0,7,8\},
				\\ \codim(\dom(T)) < \infty, \mbox{ for } i\in \{1,3,5\},
				\\ \dom(T) \mbox{ closed}, \mbox{ for } i \in \{4,6\},
				\\ \dom(T) \mbox{ complemented}, \mbox{ for } i =2,
				\\ \dom(T^n) = \dom(T^{n+1}) \mbox{ for some } n \in \NN, \mbox{ for } i=9,
			\end{cases}
			\\\sigma_i(T) \cup \{\infty\}, \mbox{ otherwise}.
		\end{cases}
	\end{align*}
\end{definition}
Note that $\widetilde{\sigma}_0(T)$ is the usual extended spectrum $\widetilde{\sigma}(T)$. If the resolvent set $\rho(T)$ is not empty, $\widetilde{\sigma}_i(T)$ coincides with the extended essential spectrum introduced by González and Onieva \cite{gonzalez1985spectral}, which satisfies that $\infty \in \widetilde{\sigma}_i(T)$ if and only if $0\in \sigma_i((\mu - T)^{-1})$ for any $\mu \in \rho(T)$. In particular, if $T$ has non-empty resolvent set, $\widetilde{\sigma}_i(T)$ are non-empty compact subsets of $\CC_\infty$ except for $i =9$ (see \cite{gramsch1971spectral}), where $\CC_\infty$ denotes the Riemann sphere $\CC \cup \{\infty\}$. 
If $T$ has empty resolvent set, $\sigma_i(T)$ is a closed subset of $\CC$ for $i\in \{0,1,2,4,5,6,7\}$, see \cite[Section I.3]{edmunds1987spectral} and \cite{yood1951properties}. We do not know if $\widetilde{\sigma_i}(T)$ or $\sigma_i(T)$ are closed in the other cases.

Now we turn to the definition of the regularized functional calculus of bisectorial-like operators. Its construction is completely analogous to the one of the regularized functional calculus of sectorial operators given by Haase in \cite{haase2005general, haase2006functional}, and the adaptation of it from the sectorial operators to the bisectorial-like operators is straightforward.

Given any $\varphi\in (0,\pi)$, we denote the sector $S_{\varphi}:=\left\{z\in\CC: \left|\mbox{arg}(z)\right|<\varphi\right\}$. For any $\omega \in (0,\pi/2]$ and $a\geq 0$,  we set the bisector
\begin{equation*}
	BS_{\omega,a} :=
	\begin{cases}
		(-a + S_{\pi-\omega}) \cap (a-S_{\pi-\omega})\;\;&\mbox{ if } \omega< \pi/2  \mbox{ or }  a> 0,\\
		i\RR &\mbox{ if } \omega =\pi/2 \mbox{ and } a=0.
	\end{cases}
\end{equation*}


\begin{definition}
	Let $(\omega,a) \in (0,\pi/2] \times [0,\infty)$ and let $A \in C(X)$.  We will say that $A$ is a {\bf bisectorial-like operator} of angle $\omega$ and half-width $a$ if the following conditions hold:
	\begin{itemize}
		\item $\sigma (A) \subseteq \overline{BS_{\omega,a}}$.
		\item For all $\omega' \in (0, \omega)$, $A$ satisfies the resolvent bound
		\begin{align*}
			\sup \Big\{ \min\{|\lambda-a|, |\lambda+a|\} \| (\lambda-A)^{-1}\| \, : \, \lambda \notin \overline{BS_{\omega',a}} \Big\}< \infty.
		\end{align*}
	\end{itemize}
	We also set $M_A := \widetilde{\sigma}(A) \cap \{-a,a,\infty\}$. For the rest of the paper, $(\omega,a)$ will denote a pair in $(0,\pi/2] \times [0,\infty)$.
\end{definition}

\noindent Given a Banach space $X$, we will denote the set of all bisectorial-like operators on $X$ of angle $\omega$ and half-width $a$ in $X$ by $\BSect(\omega,a)$. We will omit an explicit mention to $X$ for the sake of simplicity. Notice that $A \in \BSect(\omega,a)$ if and only if both $a+A, a-A$ are sectorial of angle $\pi-\omega$ in the sense of \cite{haase2006functional}.

We denote by $\mathcal{O}(\Omega), \mathcal{M}(\Omega)$ the sets of holomorphic functions and meromorphic functions defined in an open subset $\Omega \subseteq \CC$, respectively. For $A \in \BSect(\omega,a)$, let $U_A := \{-a,a,\infty\} \setminus \widetilde \sigma(A)$. If $\sigma(A) \neq \emptyset$, set
\begin{align*}
	r_d := \begin{cases}
		\mbox{dist}\{d, \sigma(A)\}, \quad & \mbox{if } d\in \{-a,a\},
		\\ r(A)^{-1}, \quad & \mbox{if } d = \infty,
	\end{cases}
	\qquad d \in U_A,
\end{align*}
where $\mbox{dist}\{\cdot, \cdot\}$ denotes the distance between two sets, and $r(A)$ the spectral radius of $A$. If $\sigma(A) = \emptyset$ (so $\widetilde{\sigma}(A) = \{\infty\}$ and $\infty \notin U_A$), set $r_a = r_{-a} := \infty$. 

For $d \in U_A$ suppose that $s_d \in (0,r_d)$. Then, for $\varphi\in (0,\omega)$, set $\Omega(\varphi, (s_d)_{d\in U_A})$ as follows. If $U_A = \emptyset$ (i.e., $M_A = \{-a,a,\infty\}$), we set $\Omega_\varphi := BS_{\varphi,a}$. Otherwise, for each $d \in U_A$, let $B_d(s_d)$ be a ball centred at $d$ of radius $s_d$, where $B_\infty(r_\infty) = \{z \in \CC \, | \, |z|>r_\infty ^{-1}\}$. Then, we set $\Omega(\varphi, (s_d)_{d\in U_A}) := BS_{\varphi,a} \setminus (\bigcup_{d\in U_A} \overline{B_d (s_d)})$. Note that, if $\varphi < \varphi' < \omega$ and $s_d < s_d' < r_d$ for each $d \in U_A$, then the inclusion $\Omega(\varphi', (s_d')_{d\in U_A}) \subseteq \Omega(\varphi, (s_d)_{d\in U_A})$ holds. Thus we can form the inductive limits
\begin{align*}
	\mathcal{O}[\Omega_A] &:= \bigcup \left\{ \mathcal{O}(\Omega(\varphi, (s_d)_{d\in U_A})) \, \Big| \, 0 < \varphi < \omega, \, 0 < s_d < r_d \mbox{ for } d \in U_A\right\}, 
	\\ \mathcal{M}[\Omega_A] &:= \bigcup \left\{ \mathcal{M}(\Omega(\varphi, (s_d)_{d\in U_A})) \, \Big| \, 0 < \varphi < \omega, \, 0 < s_d < r_d \mbox{ for } d \in U_A\right\}. 
\end{align*}

Hence, $\mathcal{O}[\Omega_A], \mathcal{M}[\Omega_A]$ are algebras of holomorphic functions and meromorphic functions (respectively) defined on an open set containing $\widetilde{\sigma}(A) \setminus M_A$. Next, we define the following notion of regularity at $M_A$.
\begin{definition}\label{regularityDef}
	Let $f \in \mathcal{M}[\Omega_A]$. We say that $f$ is regular at $d\in \{-a,a\}\cap M_A$ if $\lim_{z \to d} f(z) =: c_d \in \CC$ exists and, for some $\varphi \in (0,\omega)$
	\begin{align*}
		\int_{\partial BS_{\varphi',a}, |z-d|< \varepsilon}  \left|\frac{f(z)-c_d}{z-d}\right| |dz| < \infty, \quad \text{for some } \varepsilon > 0 \mbox{ and for all } \varphi' \in \left(\varphi, \frac{\pi}{2}\right].
	\end{align*}
	If $\infty \in M_A$, we say that $f$ is regular at $\infty$ if $\lim_{z\to \infty} f(z) =: c_\infty \in \CC$ exists and
	\begin{align*}
		\int_{\partial BS_{\varphi',a}, |z| > R}  \left|\frac{f(z)-c_\infty}{z}\right| |dz| < \infty, \quad \text{for some } R > 0 \mbox{ and for all } \varphi' \in \left(\varphi, \frac{\pi}{2}\right].
	\end{align*}
	We say that $f$ is quasi-regular at $d\in M_A$ if $f$ or $1/f$ is regular at $d$. Finally, we say that $f$ is (quasi-)regular at $M_A$ if $f$ is (quasi-)regular at each point of $M_A$.
\end{definition}

\begin{figure}[h]
	\centering
	\includegraphics[width=0.45\textwidth]{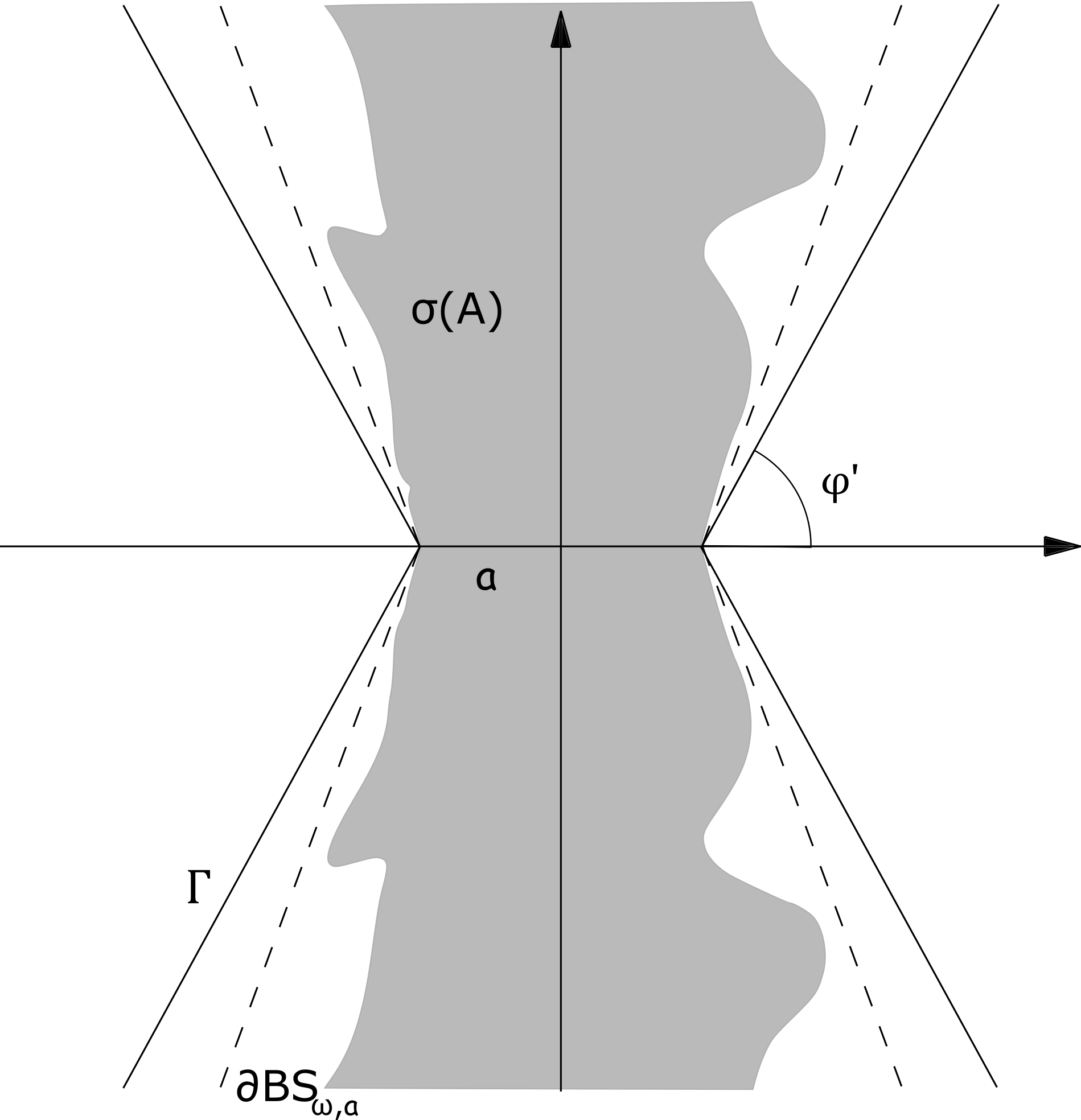}
	\caption{Spectrum of a bisectorial-like operator and integration path of the functional calculus.}
	\label{BisectFigure}
\end{figure}

\begin{remark}
	Note that if $f$ is regular at $M_A$ with every limit being not equal to $0$, then $1/f$ is also regular at $M_A$. If $f$ is quasi-regular at $M_A$, then $\mu-f$ and $1/f$ are also quasi-regular at $M_A$ for each $\mu \in \CC$. A function $f$ which is quasi-regular at $M_A$ has well-defined limits in $\CC_\infty$ as $z$ tends to each point of $M_A$.
\end{remark}


Next, let $\mathcal{E}(A)$ be the subset of functions of $\mathcal{O}[\Omega_A]$ which are regular at $M_A$. Then, for any $b \in \CC \backslash \overline{BS_{\varphi,a}}$, the set equality
\begin{align*}
	\mathcal{E}(A) &= \mathcal{E}_0(A) + \CC \frac{1}{b+z} + \CC \frac{1}{b-z} + \CC \mathbf{1},
\end{align*}
holds true, where $\mathbf{1}$ is the constant function with value $1$, and
\begin{align*}
	\mathcal{E}_0(A):= &\Bigg\{f \in \mathcal{O}[\Omega_A] \, : \, f \mbox{ is regular at } M_A \mbox{ with } \lim_{z \to d} f(z) = 0 \mbox{ for all } d \in M_A \Bigg\}.
\end{align*}

Given a bisectorial-like operator $A \in \BSect (\omega,a)$, we define the algebraic homomorphism $\Phi: \mathcal{E}(A) \to \mathcal{L}(X)$ determined by setting $\Phi\left(\frac{1}{b+z}\right) = (b+A)^{-1}$, $\Phi\left(\frac{1}{b-z}\right) = (b-A)^{-1}$, $\Phi(\mathbf{1}) = I$, and
\begin{align*}
	\Phi(f) := f(A) := \frac{1}{2\pi i} \int_\Gamma f(z) (z-A)^{-1}\, dz, \quad f \in \mathcal{E}_0(A),
\end{align*}
where $\Gamma$ is the positively oriented boundary of $\Omega(\varphi', (s_d')_{d\in U_A})$ with $\varphi < \varphi' < \omega$ and $s_d < s_d' < r_d$, where $f \in \mathcal{O}(\Omega(\varphi, (s_d)_{d\in U_A}))$. Note that the above integral is well defined in the Bochner sense since $z \mapsto (z-A)^{-1}$ is analytic, so continuous, and  $\int_\Gamma |f(z)| \|(z-A)^{-1}\|\, |dz| < \infty$. It is readily seen that $f(A)$ is well defined for $f \in \mathcal{E}(A)$ (that is, that $f(A)$ does not depend on the election of $\varphi', (s_d')_{d\in U_A}$).

Next, we follow the regularization method given in \cite{haase2005general} to extend the functional calculus $\Phi$ to a regularized functional calculus (also denoted by $\Phi$), which involves meromorphic functions.
\begin{definition}
	Let $a\geq 0, \, 0< \omega \leq \frac{\pi}{2}$ and $A \in \BSect(\omega,a)$. Then, a function $f \in \mathcal{M}[\Omega_A]$ is called regularizable by $\mathcal{E} (A)$ if there exists $e \in \mathcal{E}(A)$ such that
	\begin{itemize}
		\item $e(A)$ is injective,
		\item $ef \in \mathcal{E}(A)$.
	\end{itemize}
	For any regularizable $f \in \mathcal{M}[\Omega_A]$ with regularizer $e \in \mathcal{E}(A)$, we set
	\begin{align*}
		\Phi(f) := f(A):= e(A)^{-1} (ef)(A).
	\end{align*}
\end{definition}
By \cite[Lemma 3.2]{haase2005general}, one has that this definition is independent of the regularizer $e$, and that $f(A)$ is a well-defined closed operator. We will denote by $\mathcal{M}(A)$ the subset of functions of $\mathcal{M}[\Omega_A]$ which are regularizable by $\mathcal{E}(A)$. As in the case for sectorial operators \cite[Theorem 3.6]{haase2005general}, this regularized functional calculus satisfies the properties given in the lemma below. Recall that, for $A,B \in C(X)$, we mean by $A \subseteq B$ that $\mathcal D(A) \subseteq \mathcal D(B)$ with $Ax = Bx$ for every $x \in \mathcal D(A)$.

\begin{lemma}\label{functionalCalculusProperties}
	Let $A \in \BSect(\omega,a)$ and $f \in \mathcal{M}(A)$. Then
	\begin{enumerate}
		\item If $T\in \mathcal{L}(X)$ commutes with $A$, that is, $TA \subseteq AT$, then $T$ also commutes with $f(A)$, i.e. $T f(A) \subseteq f(A)T$.
		\item $\zeta(A) = A$, where $\zeta(z) = z$, $z \in \CC$.
		\item Let $g \in \mathcal{M}(A)$. Then
		\begin{align*}
			f(A) + g(A) \subseteq (f+g)(A), \qquad f(A) g(A) \subseteq (f g)(A).
		\end{align*}
		Furthermore, $\dom(f(A)g(A)) = \dom((fg)(A)) \cap \dom(g(A))$, and one has equality in these relations if $g(A) \in \mathcal{L}(X)$.
		\item Let $\lambda \in \CC$. Then 
		\begin{align*}
			\frac{1}{\lambda - f(z)} \in \mathcal{M}(A) \iff \lambda - f(A) \text{ is injective}.
		\end{align*}
		If this is the case, $(\lambda-f(z))^{-1}(A) = (\lambda -f(A))^{-1}$. In particular, $\lambda \in \rho(A)$ if and only if $(\lambda - f(z))^{-1} \in \mathcal{M}(A)$ with $(\lambda - f(A)^{-1}) \in \mathcal{L}(X)$.
	\end{enumerate}
\end{lemma}
\begin{proof}
	The statement follows by straightforward applications of the Cauchy's theorem, the resolvent identity, and \cite[Section 3]{haase2005general}.
\end{proof}

\begin{lemma}\label{pointSpectraLemma}
	Let $A \in \BSect(\omega,a)$ and $f \in \mathcal{M}(A)$. Then $f(A)x = f(\lambda)x$ for any $x \in \kernel (\lambda - A)$.
\end{lemma}
\begin{proof}
	See \cite[Proposition 3.1]{haase2005spectral} for the analogous result for sectorial operators.
\end{proof}

\begin{lemma}\label{regularizerLemma}
	Let $A \in \BSect(\omega,a)$, $f \in \mathcal{M}(A)$ and $\lambda \in \widetilde{\sigma}(A) \setminus M_A$ such that $f(\lambda) \neq \infty$. There is a regularizer $e \in \mathcal{E}(A)$ for $f$ with $e(\lambda) \neq 0$.
\end{lemma}

\begin{proof}
	The proof is analogous to the case of sectorial operators, see \cite[Lemma 4.3]{haase2005spectral}.
\end{proof}

\begin{lemma}\label{differentNFCLemma}
	Let $A\in \BSect(\omega,a)$ and $f \in \mathcal{M}[\Omega_A]$. Assume that $f$ is regular at $M_A$ and that all the poles of $f$ are contained in $\CC \backslash \sigma_p(A)$. Then, $f \in \mathcal{M}(A)$.
	Moreover, if every pole of $f$ is contained in $\rho(A)$, then $f(A)\in \mathcal{L}(X)$.
\end{lemma}

\begin{proof}
	The proof is the same as in the case of sectorial operators, see \cite[Lemma 6.2]{haase2005spectral}. We include it here since we need it in the proof of Theorem \ref{boundedFunctionalCTheorem}.
	
	Let $f \in \mathcal{M}[\Omega_A]$ be as required. That is, there exists $\varphi \in (0,\omega)$ and $s_d \in (0,r_d)$ for each $d \in U_A$ such that $f \in \mathcal{M}(\Omega(\varphi, (s_d)_{d\in U_A}))$. 
	Since $f$ has finite limits at $M_A$, we can assume that $f$ has only finitely many poles by making $\varphi, (s_d)_{d\in U_A}$ bigger. 
	Thus, let $\lambda_j$ for $j\in \{1,...,N\}$ be an enumeration of those poles of $f$ and let $n_j \in \NN$ be the order of pole of $f$ located at $\lambda_j$, for $j \in \{1,...,N\}$. Then, the function $\displaystyle{g(z) := f(z) \prod_{j=1}^{N}\frac{(\lambda_j-z)^{n_j}}{(b-z)^{n_j}}}$ has no poles, i.e. $g\in \mathcal{O}[\Omega_A]$, and is regular at $M_A$. Hence $g\in \mathcal{E}(A)$. Moreover, setting $\displaystyle{r(z) := \prod_{j=1}^{N}\frac{(\lambda_j-z)^{n_j}}{(b-z)^{n_j}}}$, one has that $\displaystyle{r(A) = \prod_{j=1}^{N}(\lambda_j-A)^{n_j}(b-A)^{-n_j}}$ is bounded and injective, since by assumption $\{\lambda_1,...,\lambda_n\} \subseteq \CC \setminus \sigma_p(A)$. In short, $f$ is regularized by $r$, so $f \in \mathcal{M}(A)$.
	
	Now, assume that the poles of $f$ lie inside $\rho(A)$. Then the operator $r(A)$ is not only bounded and injective, but invertible too, from which follows that $f(A) = r(A)^{-1} (rf)(A) \in \mathcal{L}(X)$.
\end{proof}

\section{Spectral mapping theorems for $M_A = \emptyset$}\label{noSingularPointsSection}


For $A \in \BSect(\omega,a)$, the spectral mapping theorems \eqref{spectralMappingTheoremIntro} given in \cite{gonzalez1985spectral, gramsch1971spectral} are applicable to every $f \in \mathcal{E}(A)$ whenever $M_A = \emptyset$. This section is devoted to extend these spectral mapping theorems for all $f \in \mathcal{M}(A)$.

First, we proceed to state the spectral mapping inclusion of the spectrum $\widetilde{\sigma}$.

\begin{proposition}\label{spectralInclusionProp}
	Let $A \in \BSect(\omega,a)$, $f \in \mathcal{M}(A)$, and assume that $f$ is quasi-regular at $M_A$. Then 
	\begin{align*}
		\widetilde{\sigma}(f(A)) \subseteq f(\widetilde{\sigma}(A)).
	\end{align*}
\end{proposition}
\begin{proof}
	The proof runs along the same lines as in the case of sectorial operators, see \cite[Proposition 6.3]{haase2005spectral}. As in Lemma \ref{differentNFCLemma}, we include the proof here since it will be needed in the proof of Theorem \ref{boundedFunctionalCTheorem}.
	
	Take $\mu \in \CC$ such that $\mu \notin f(\widetilde{\sigma}(A))$. Then $\displaystyle{\frac{1}{\mu - f} \in \mathcal{M}[\Omega_A]}$ is regular in $M_A$, and all of its poles are contained in $\rho(A)$. By Lemma \ref{differentNFCLemma}, we conclude that $(\mu - f)^{-1} \in \mathcal{M}(A)$ with $(\mu - f)^{-1}(A)$ is a bounded operator. Thus, it follows that $\mu - f(A)$ is invertible, hence $\mu \notin \widetilde{\sigma}(f(A))$.
	
	Assume now that $\mu = \infty \notin f(\widetilde{\sigma}(A))$. Then $f$ is regular at $M_A$ and its poles are contained in $\rho(A)$. Another application of Lemma \ref{differentNFCLemma} yields that $f(A)$ is a bounded operator, so $\infty \notin \widetilde{\sigma}(f(A))$.
\end{proof}

Next, we give some technical lemmas.
\begin{lemma}\label{rangeLemma}
	Let $A \in \BSect(\omega,a)$, $e\in \mathcal{M}(A)$ with $e(A) \in \mathcal{L}(X)$ injective, $\lambda, b \in \CC$ with $b \in \rho(A)$. Assume that there is $c\in \CC\setminus\{0\}$ such that
	$$f(z):= \frac{b-z}{\lambda-z} (e(z)-c) \in \mathcal{M}(A) \qquad \mbox{ with } f(A) \in \mathcal{L}(X).
	$$
	Then $\mathcal{R}(\lambda- A) = \mathcal{R}\left((\lambda - A)(b-A)^{-1} e(A)^{-1}\right) = \mathcal{R}\left(e(A)^{-1} (\lambda - A)(b-A)^{-1}\right)$.
\end{lemma}

\begin{proof}
	Note that $\displaystyle{\mathcal{R}(\lambda- A) = \mathcal{R}\left((\lambda - A)(b-A)^{-1}\right) = \mathcal{R}\left((\lambda - A)(b-A)^{-1} e(A)^{-1}\right)}$ since $e(A)^{-1}$ is surjective. Moreover, since $\displaystyle{(\lambda - A)(b-A)^{-1} e(A)^{-1} \subseteq e(A)^{-1}(\lambda - A)(b-A)^{-1}}$ by Lemma \ref{functionalCalculusProperties}, we have $\displaystyle{\mathcal{R}\left((\lambda - A)(b-A)^{-1} e(A)^{-1}\right) \subseteq \mathcal{R}\left(e(A)^{-1} (\lambda - A)(b-A)^{-1}\right)}$. Thus, all is left to prove is the reverse inclusion.
	
	Let $\displaystyle{u \in \mathcal{R} \left(e(A)^{-1} (\lambda - A)(b-A)^{-1}\right)}$, so there is $x \in X$ such that $\displaystyle{e(A)u = (\lambda - A)(b-A)^{-1} x}$. Since $\displaystyle{e(z) = \frac{\lambda-z}{b-z}f(z) +c}$, one has that $\displaystyle{u = \frac{1}{c}(\lambda - A)(b-A)^{-1}(x - f(A)u)}$ so $u \in \ran(\lambda - A) = \displaystyle{\ran\left((\lambda - A)(b-A)^{-1} e(A)^{-1}\right)}$, and the claim follows.
\end{proof}

\begin{lemma}\label{domainLemma}
	Let $f \in \mathcal{M}(A)$ and $\lambda \in \sigma(A)\setminus M_A$ with $f(\lambda) = 0$, and let $b \in \rho(A)$. If $\displaystyle{g(z) := \frac{b-z}{\lambda-z} f(z)}$, then $g\in \mathcal{M}(A)$, $\dom(g(A)) = \dom(f(A))$ and
	$$f(A) = (\lambda - A)(b-A)^{-1} g(A) = g(A) (\lambda - A)(b-A)^{-1}.
	$$
\end{lemma}
\begin{proof}
	Let $e \in \mathcal{E}(A)$ be a regularizer for $f$ with $e(\lambda)\neq 0$, see Lemma \ref{regularizerLemma}. The fact that $eg$ has the same behaviour as $ef$ at $M_A$ implies that $eg \in \mathcal{E}(A)$, that is, $e$ is a regularizer for $g$ and $g \in \mathcal{M}(A)$, so $g(A)$ is well defined.
	
	On one hand, it follows by Lemma \ref{functionalCalculusProperties} (3) that $f(A) = g(A) (\lambda - A)(b-A)^{-1} \supset (\lambda - A)(b-A)^{-1} g(A)$ with $\dom\left((\lambda - A)(b-A)^{-1} g(A)\right) = \dom(f(A)) \cap \dom(g(A))$. By the definition of composition of closed operators, $\dom\left((\lambda - A)(b-A)^{-1} g(A)\right) = \dom(g(A)) \cap g^{-1}\left(\dom\left((\lambda - A)(b-A)^{-1}\right)\right) = \dom(g(A))$ since $\dom\left((\lambda - A)(b-A)^{-1}\right) = X$. As a consequence, $\dom(g(A)) \subseteq \dom(f(A))$.
	
	Next, let $x \in \dom(f(A))$ and set $\widetilde{x}:=(eg)(A)x$. One has
	\begin{align*}
		f(A) &= e(A)^{-1} \left(\frac{\lambda-z}{b-z} (eg)(z)\right)(A) = e(A)^{-1} (\lambda - A)(b-A)^{-1} (eg)(A),
	\end{align*}
	so $\displaystyle{\widetilde{x} \in \dom\left(e(A)^{-1} (\lambda - A)(b-A)^{-1}\right)}$. An application of Lemma \ref{rangeLemma} with $c = e(\lambda) \neq 0$ shows that there is $v \in \mathcal{R}(e(A))$ with $\displaystyle{(\lambda - A)(b-A)^{-1} e(A)^{-1} v = e(A)^{-1} (\lambda - A)(b-A)^{-1} \widetilde{x}}$. By composing with $e(A)$ one gets $\widetilde{x} - v \in \kernel \displaystyle{\left((\lambda - A)(b-A)^{-1}\right)} = \kernel (\lambda - A)$. Moreover, $\kernel (\lambda - A) \subseteq \mathcal{R}(e(A))$ since $y = \frac{1}{e(\lambda)} e(A)y$ for any $y \in \kernel (\lambda-A)$, see Lemma \ref{pointSpectraLemma}. Hence, $\widetilde{x} = (\widetilde{x}-v)+v \in \mathcal{R}(e(A))$, that is $x \in \dom(g(A))$, so $\dom(g(A)) = \dom(f(A))$. Lemma \ref{differentNFCLemma}(3) shows that
	$$f(A) =  (\lambda - A)(b-A)^{-1} g(A) = g(A) (\lambda - A)(b-A)^{-1},
	$$
	and the commutativity property follows.
\end{proof}

\begin{remark}\label{i89EquivalenceRemark}
	Let $T \in C(X)$ with non-empty resolvent set, and $\ascent(T), \descent(T) < \infty$. Then $\ascent(T) = \descent(T)=:p_T$ and $X = \kernel(T^{p_T}) \oplus \ran(T^{p_T})$, see for example \cite[Theorem V.6.2]{taylor1958introduction}.
\end{remark}

Lemma below is inspired by \cite[Lemma 5]{gramsch1971spectral}.
\begin{lemma}\label{i89Lemma}
	Let $S,T \in C(X)$ with non-empty resolvent set and such that $ST = TS$. One has the following
	\begin{itemize}
		\item [(a)] If $S,T \in \Phi_9$, then $ST \in \Phi_9$. If $\ran(T) \subset \dom(S)$ and $S,T \in \Phi_8$, then $ST \in \Phi_8$.
		\item [(b)] Assume that $T$ is injective and $\dom(S) \subset \ran(T)$. If $ST \in \Phi_i$, then $S \in \Phi_i$ for $i \in \{8,9\}$.
	\end{itemize}
\end{lemma}
\begin{proof}
	(a) Let $S,T \in \Phi_9$. By Remark \ref{i89EquivalenceRemark}, $\ascent(S) = \descent(S) =: p_S$, $\ascent(T) = \descent(T) =: p_T$, and $X = \ran(S^{p_S}) \oplus \kernel(S^{p_S}) = \ran(T^{p_T}) \oplus \kernel(T^{p_T})$. Let $P_S$ be the projection onto $\kernel(S^{p_S})$ along $\ran(S^{p_S})$, $Q_S := I - P_S$, and set the analogous projections $P_T, Q_T$. Since $ST = TS$, one has that $\kernel (T^n) \subset \dom(S)$, $\kernel(S^n) \subset \dom(T)$ for any $n \in \NN$, and that $P_S, Q_S, P_T, Q_T$ commute between themselves. Then $Q := Q_T Q_S$ is a bounded projection onto $\ran(S^{p_S}) \cap \ran(T^{p_T})$, and it is readily seen that $ST$ is a (possibly unbounded) invertible operator when restricted to $Q(X)$. Since $I-Q = P_S + P_T - P_S P_T$, it is clear that $(I-Q)(X) \subset \kernel (ST)^{\max\{p_S, p_T\}}$. Then, $ST \in \Phi_9$ with $\ascent(ST) = \descent(ST) \leq \max\{p_S, p_T\} < \infty$ by \cite[Problem V.6]{taylor1958introduction}. 
	
	If in addition, $S,T \in \Phi_8 \subset \Phi_7$ with $\ran(T) \subset \dom(S)$, then $ST \in \Phi_7$, see \cite[Theorem I.3.16]{edmunds1987spectral} (although this result is stated for bounded operators, its proof is purely algebraic). Hence, we conclude that $ST \in \Phi_8$.
	
	(b) It follows by induction that $\dom(S^n) \subset \ran(T^n)$ for $n \in \NN$. Since $(ST)^n = S^n T^n = T^n S^n$ and $T$ is injective, one has that $\ascent(ST) = \ascent(S)$ and $\descent(ST) = \descent(S)$, and the claim follows (note that $ST \in \Phi_1$ implies that $S,T \in \Phi_1$). 
\end{proof}

\begin{remark}\label{resolventRemark}
	Let $T \in C(X)$ with non-empty resolvent set, and let $b \in \rho(T)$ and $\lambda \in \CC$. Then $\lambda - T \in \Phi_i$ if and only if $\displaystyle{(\lambda - T)(b-T)^{-1}} \in \Phi_i$ for $i \in \{0,1,2,3,4,5,6,7,8,9\}$, see for example \cite[Lemma 1]{gonzalez1985spectral}.
\end{remark}

\begin{lemma}\label{fEquivalentgLemma}
	Let $A\in \BSect(\omega,a), \, f,g \in \mathcal{M}(A)$ with $f,g$ quasi-regular at $M_A$,  $0 \notin g(\widetilde{\sigma}(A))$ and such that
	$$f(z) := g(z)\prod_{j=1}^N \left(\frac{\lambda_j-z}{b-z}\right)^{n_j}, 
	$$
	for some $b \in \rho(A), \, \lambda_j \in \sigma(A)\setminus M_A,$ and $n_j \in \NN$ for $j=1,...,N$. Then 
	\begin{enumerate}
		\item[(a)] If $f(A) \in \Phi_i$, then $\lambda_j - A \in \Phi_i$ for all $j=1,...,N$ and for $i \in \{0,1,2,3,4,5,6,8,9\}$.
		\item[(b)] If $\lambda_j - A \in \Phi_i$ for all $j=1,...,N$, then $f(A) \in \Phi_i$ for $i \in \{0,1,2,3,4,5,7,8,9\}$.
	\end{enumerate}
\end{lemma}
\begin{proof}
	Set $r(z):= \displaystyle{\prod_{j=1}^N \left(\frac{\lambda_j-z}{b-z}\right)^{n_j}}$, so $r(A) \in \mathcal{L}(X)$. Several applications of Lemma \ref{domainLemma} imply that $\dom(g(A))  = \dom(f(A))$ and $f(A) = r(A) g(A)$. Moreover, Proposition \ref{spectralInclusionProp} yields that $0 \notin \widetilde{\sigma}(g(A))$, so $g(A)$ is surjective and injective. Thus $g(A): \dom (g(A)) = \dom(f(A)) \to X$ is an isomorphism when $\dom(f(A))$ is endowed with the graph norm given by $f(A)$ (which is equivalent to the graph norm given by $g(A)$). Therefore, $f(A) \in \Phi_i$ if and only if $r(A)$. This follows by the very definition of $\Phi_i$ for all $i \in \{0,1,2,3,4,5,6,7\}$, and  by Lemma \ref{i89Lemma} for $i \in \{8,9\}$. 
	
	
	Since the bounded operators $\displaystyle{(\lambda_j-A)(b-A)^{-1}}$ commute between themselves, we have:
	\begin{enumerate}
		\item if $r(A) \in \Phi_i$, then $\displaystyle{(\lambda_j - A)(b-A)^{-1} \in \Phi_i}$ for all $j=1,...,N$, and for $i \in \{0,1,2,3,4,5,6,8,9\}$,
		\item If $\displaystyle{(\lambda_j - A)(b-A)^{-1} \in \Phi_i}$ for all $j=1,...,N$, then $r(A) \in \Phi_i$, for $i \in \{0,1,2,3,4,5,7,8,9\}$.
	\end{enumerate}
	see for example \cite[Lemma 3]{gonzalez1985spectral} and \cite[Lemma 5(c)]{gramsch1971spectral}. 
	Hence, the claim follows from Remark \ref{resolventRemark}.

	
\end{proof}

We give now the main result of this section.
\begin{proposition}\label{essentialMappingNoSingular}
	Let $A \in \BSect(\omega,a)$, $f \in \mathcal{M}(A)$, where $f$ is quasi-regular at $M_A$. Then
	\begin{enumerate}
		\item[(a)] $f(\widetilde{\sigma}_i(A) ) \setminus f(M_A)\subset \widetilde{\sigma}_i(f(A))$ for $i\in \{0,1,2,3,4,5,6,8,9\}$.
		\item [(b)] $\widetilde{\sigma}_i (f(A)) \subset f(\widetilde{\sigma}_i)(A) \cup f(M_A)$ for $i\in \{0,1,2,3,4,5,7,8,9\}$.
	\end{enumerate}
	
\end{proposition}
\begin{proof}
	Take $i \in \{0,1,2,3,4,5,6,8,9\}$ and let $\mu \in \CC$ be such that $\mu \in f(\widetilde{\sigma}_i(A))\setminus f(M_A)$. By considering the function $f - \mu$ instead of $f$, we can assume without loss of generality that $\mu = 0$. As $0\notin f(M_A)$, $f^{-1}(0) \cap \widetilde{\sigma}(A)$ must be finite. Let $\lambda_1,\ldots, \lambda_N$ be all the points in $f^{-1}(0) \cap \widetilde{\sigma}(A)$ (so $\lambda_j \in \widetilde{\sigma}_i(A)$ for some $j \in \{1,...,N\}$), and let $n_j$ be the order of the zero of $f$ at $\lambda_j$. Let $b \in \rho(A)$ and set
	\begin{equation}\label{gDefinition1}
		g(z) := f(z) \prod_{j=1}^N \left(\frac{b-z}{\lambda_j-z}\right)^{n_j}.
	\end{equation}	
	Then $0 \notin g(\widetilde{\sigma}(A))$ and is $g$ quasi-regular at $M_A$. Several applications of Lemma \ref{domainLemma} imply that $g \in \mathcal{M}(A)$, and Lemma \ref{fEquivalentgLemma}(a) yields that $f(A) \notin \Phi_i$.
	
	Take now $i \in \{0,1,2,3,4,5,7,8,9\}$ and let $\mu \in \CC$ be such that $\mu \notin f(\widetilde{\sigma}_i(A)) \cup f(M_A)$. We will prove that $\mu \notin \widetilde{\sigma}_i(f(A))$. We can assume $\mu=0$. Again, $f^{-1}(0) \cap \sigma(A)$ has finite cardinal, so let $g$ be as given in \eqref{gDefinition1}. Since $\lambda_j - A \in \Phi_i$ for all $j=1,...,n$, applications of Lemma \ref{domainLemma} and Lemma \ref{fEquivalentgLemma}(b) yield that $f(A) \in \Phi_i$, as we wanted to show.

	
	Assume now that $\mu = \infty$. If $\rho(f(A)) \neq \emptyset$ take $b \in \rho(f(A))$. An application of what we have already proven to the function $\displaystyle{\frac{1}{b-f(z)}}$ shows the claim, see the paragraph below Definition \ref{extendedEssSpectraDef}. Hence, all that is left to prove is that we can assume without loss of generality that $\rho(f(A)) \neq \emptyset$. Take $\nu \in \CC\setminus f(M_A)$, so $f^{-1}(\nu) \cap \widetilde{\sigma}(A)$ has finite cardinal. Let $\nu_1,\ldots,\nu_M$ be all the points in $f^{-1}(\nu) \cap \sigma(A)$, and let $m_j$ be the order of the zero of $f-\nu$ at $\nu_j$. Let $b \in \rho(A)$ and set
	\begin{equation}\label{gdefinition}
		h(z) := (f(z)-\nu) \prod_{j=1}^M \left(\frac{b-z}{\nu_j-z}\right)^{m_j}.
	\end{equation}
	Lemma \ref{domainLemma} yields that $h \in \mathcal{M}(A)$ with $\dom(f(A)) = \dom(h(A))$, and using \eqref{gdefinition} it is readily seen that $\dom(f(A)^n) = \dom(h(A)^n)$ for all $n \in \NN$. In particular, $\infty \in \widetilde{\sigma}_i(f(A))$ if and only if $\infty \in \widetilde{\sigma}_i(h(A))$. Since $0 \notin h(\widetilde{\sigma}(A))$, Proposition \ref{spectralInclusionProp} implies that $0 \in \rho(h(A))$. Therefore, we can assume that $\rho(f(A)) \neq \emptyset$, and the proof is done.
\end{proof}



\section{General case}\label{spectralMappingSection}
In this section we deal with the case $M_A \neq \emptyset$. The difficulty of this setting arises from the fact that $f$ is not necessarily either holomorphic or meromorphic at $M_A$, so the factorization techniques used in Section \ref{noSingularPointsSection} do not apply here.

First, we give some remarks about $M_A$ which will be the key for the proof of the spectral mapping theorems.
\begin{remark}\label{isolatedEssentialPoints}
	Let $T\in C(X)$ with non-empty resolvent set, $d \in \widetilde{\sigma}(T)$ with $d$ an accumulation point of $\rho(T)$, and $i \in \{1,2,3,4,5,6,7,8\}$. The following statements about the essential spectrum are well known, see for example \cite[Sections I.3 \& I.4]{edmunds1987spectral}, \cite[Chapter 4\S5]{kato1966perturbation} and \cite[Section V.6]{taylor1958introduction}.
	\begin{itemize}
		\item[(a)] If $d$ is also an accumulation point of $\widetilde{\sigma}(T)$, then $d \in \widetilde{\sigma}_i(T)$.
		\item[(b)] If $d \in \widetilde{\sigma}_i(T)$ and $d$ is not an accumulation point of $\widetilde{\sigma}_i(T)$, then there is a neighborhood $\Omega$ of $d$ such that $\widetilde \sigma(T) \cap \Omega$ consists of $d$ and a countable (possibly empty) set of eigenvalues of $T$ with finite dimensional eigenspace, which are isolated between themselves.
		\item[(c)] If $d \notin \widetilde{\sigma}_i(T)$, then $d$ is an isolated point of $\widetilde{\sigma}(T)$. Moreover, $d \in \sigma_p(T)$ with $\nul(d-T) = \defec(d-T) < \infty$, $\ascent(d-T) = \descent(d-T) < \infty$, and $\dim\left(\cup_{n\geq 1}\kernel((d-T)^n)\right)< \infty$.
	\end{itemize}
\end{remark}

\begin{lemma}\label{singularPointsEquivalenceLemma}
	Let $A \in \BSect(\omega,a)$, $d \in M_A$ and $i,j \in \{1,2,3,4,5,6,7,8\}$.  Then
	\begin{itemize}
		\item $d \in \widetilde{\sigma}_i(A)$ if and only if $d \in \widetilde{\sigma}_j(A)$,
		\item if $\infty \in \widetilde{\sigma}(A)$, then $\infty \in \widetilde{\sigma}_i(A)$.
	\end{itemize} 
\end{lemma}
\begin{proof}
	If $d \in \widetilde{\sigma}_6(A)$, then $d \in \widetilde{\sigma}_i(A)$ since $\widetilde{\sigma}_6(A) \subseteq \widetilde{\sigma}_i(A)$ for any $i \in \{1,2,3,4,5,6,7,8\}$. If $d \notin \widetilde{\sigma}_6(A)$, then Remark \ref{isolatedEssentialPoints}(c) implies that $d \notin \widetilde{\sigma}_i(A)$ for $i \in \{1,2,3,4,5,6,7,8\}$, and the claim follows.
\end{proof}

Take $T \in C(X)$ with non-empty resolvent set, and let $\Lambda_1,\ldots,\Lambda_n$ be the components of $\widetilde \sigma(T)$. Let $\Lambda$ be subset of $\widetilde \sigma(T)$ which is open and closed in the relative topology of $\widetilde \sigma(T)$ (i.e. $\Lambda$ is the union of some components of $\widetilde \sigma(T))$. If $\infty \notin \Lambda$, the spectral projection $P_\Lambda$ of $T$ is given by
\begin{equation}\label{spectralProjectEq}
	P_\Lambda := \int_\Gamma (z-A)^{-1}\, dz,
\end{equation}
where $\Gamma$ is a finite collection of paths contained in $\rho(T)$ such that $\Gamma$ has index $1$ with respect to every point in $\Lambda$, and has index $0$ with respect to every point in $\sigma(T) \setminus \Lambda$. If $\infty \in \Lambda$, then the spectral projection $P_\Lambda$ of $T$ is given by $P_\Lambda := I - P_{\widetilde \sigma(T) \setminus \Lambda}$, where $P_{\widetilde \sigma(T)\setminus \Lambda}$ is as in \eqref{spectralProjectEq}.

We collect in the form of a lemma some well-known results about spectral projections, see for instance \cite[Section V.9]{dunford1963linear}.

\begin{lemma}\label{commutingProjectionsLemma}
	Let $T, \Lambda$ be as above. Then
	\begin{enumerate}
		\item $P_\Lambda$ is a bounded projection commuting with $T$;
		\item $\widetilde \sigma(T_\Lambda) = \Lambda$, where $T_\Lambda :\ran(P_\Lambda) \to \ran(P_\Lambda)$ is the part of $T$ in $\ran(P_\Lambda)$.
	\end{enumerate}
	As a consequence, for $\lambda \in \Lambda \cap \CC$, $\kernel (\lambda - T) \subseteq \ran(P_\Lambda)$ and $\ran(I - P_\Lambda) \subseteq \ran(\lambda - T)$. Also, if $\infty \notin \Lambda$, then $\ran(P_\Lambda) \subseteq \dom(T)$.
\end{lemma}

We also need the following two lemmas. \textcolor{red}{For a Jordan curve $\Gamma$, let $\mbox{\normalfont{ext}}(\Gamma),\, \mbox{\normalfont{int}}(\Gamma)$ denote the exterior and the interior of $\Gamma$ respectively.}

\begin{lemma}\label{resolventLemma}
	Let $A \in \BSect(\omega,a)$, and let $\Lambda\subset \widetilde \sigma(A)$ be an open and closed subset in the relative topology of $\widetilde \sigma(A)$. Then $A_\Lambda \in \BSect(\omega,a)$ with $\mathcal{M}(A) \subseteq \mathcal{M}(A_\Lambda)$, and one has
	$$f(A_\Lambda) = f(A)|_{\ran(P_\Lambda)}, \qquad \mbox{and} \qquad \widetilde{\sigma}_i(f(A_\Lambda)) \subseteq \widetilde{\sigma}_i(f(A)),
	$$
	for every $f \in \mathcal{M}(A)$ and $i \in \{0,1,2,3,4,5,6,9\}$.
\end{lemma}
\begin{proof}
	It follows by Lemma \ref{commutingProjectionsLemma} that $\widetilde\sigma(A_\Lambda) = \Lambda \subseteq \overline{BS(\omega,a)} \cup \{\infty\}$. Moreover, it is readily seen that $(z-A_\Lambda)^{-1} = (z-A)^{-1}|_{\ran(P_\Lambda)}$ for every $z \in \rho(A)$. As a consequence, one gets that $A_\Lambda$ is indeed a bisectorial-like operator on $\ran(P_\Lambda)$ of angle $\omega$ and half-width $a$, and that $\mathcal E(A) \subseteq \mathcal E (A_\Lambda)$ with $f(A)|_{\ran(P_\Lambda)} = f(A_\Lambda)$ for all $f \in \mathcal{E}(A)$. Thus, if $e \in \mathcal{E}(A)$ is a regularizer for $f \in \mathcal{M}(A)$, then $e$ is also a regularizer for $f$ with respect to $A_\Lambda$, so $\mathcal{M}(A) \subseteq \mathcal{M}(A_\Lambda)$.
	
	Now, we have $P_\Lambda f(A) \subseteq f(A)P_\Lambda$ for every $f \in \mathcal M(A)$ by Lemma \ref{functionalCalculusProperties}. From this and the above properties, it is not difficult to get, for every $f \in \mathcal{M}(A)$, $f(A_\Lambda) = f(A)|_{\ran(P_\Lambda)}$, with $\dom(f(A_\Lambda)) = \dom(f(A))\cap \ran(P_\Lambda), \, \kernel(f(A_\Lambda)) = \kernel(f(A)) \cap \ran(P_\Lambda), \, \ran(f(A_\Lambda)) = \ran (f(A))\cap \ran(P_\Lambda)$. Hence $\widetilde{\sigma}_i(f(A_\Lambda)) \subseteq \widetilde{\sigma}_i(f(A))$ for $i \in \{0,1,4,5,6\}$. 
	
	Thus, all that is left to prove are the spectral inclusions for $\widetilde\sigma_2, \widetilde \sigma_3$ and $\widetilde\sigma_9$. Regarding the first case, assume $f(A) \in \Phi_2$, i.e., $\ran(f(A)) \oplus W = X$ for some closed linear subspace $W$. Then $\ran(f(A_\Lambda)) \oplus (\ran(f(A)) \cap \ran(I- P_\Lambda)) \oplus W = X$ since $\ran(f(A)) = \ran(f(A_\Lambda))  \oplus  (\ran(f(A) \cap \ran(I-P_\Lambda))$. Thus, there exists a bounded projection $\tilde{P}$ from $X$ onto $\ran(f(A_\Lambda))$. Then, $\tilde P|_{\ran(P_\Lambda)}$ is a bounded projection from $\ran(P_\Lambda)$ onto $\ran(f(A_\Lambda))$, that is, $\ran(f(A_\Lambda))$ is complemented in $\ran(P_\Lambda)$, so $f(A_\Lambda) \in \Phi_2$. An analogous reasoning with $\dom(f(A))$ and $\dom(f(A_\Lambda))$ shows that, if $\infty \in \widetilde{\sigma}_2(f(A_\Lambda))$, then $\infty \in \widetilde{\sigma}_2(f(A))$. Hence, the inclusion $\widetilde{\sigma}_2(f(B)) \subseteq \widetilde{\sigma}_2(f(A))$ holds for any $f \in \mathcal{M}(A)$.
	
	Similar reasoning proves the inclusion for $\widetilde \sigma_3$. For $\widetilde\sigma_9$, the inclusion follows from $\ascent(f(A)) = \max\{\ascent(f(B)), \ascent(f(A)|_Z)\}$ and $\descent(f(A)) = \max\{\descent(f(B)), \descent(f(A)|_Z)\}$, see \cite[Problem V.6]{taylor1958introduction}.

\end{proof}

Let $A \in \BSect(\omega,a)$. Recall that $M_A \setminus \widetilde{\sigma}_i(A) =  M_A \setminus \widetilde{\sigma}_j(A)$ for $i,j \in \{1,2,...,8\}$, see Lemma \ref{singularPointsEquivalenceLemma}. Then, by Remark \ref{isolatedEssentialPoints}(c), $M_A \setminus \widetilde{\sigma}_i(A)$ is an open and closed subset of $\widetilde\sigma(A)$ (in the relative topology of $\widetilde \sigma(A)$) for $i \in \{1,2,3,4,5,6,7,8\}$. Also,  

\begin{lemma}\label{quotientNFCLemma}
	Let $A \in \BSect(\omega,a)$, and let $\Lambda = \widetilde\sigma(A) \setminus (M_A \setminus \widetilde{\sigma}_i(A))$ for any $i \in \{1,2,3,4,5,6\}$. Then $M_{A_\Lambda} \subset \widetilde{\sigma}_i(A_\Lambda)$ for $i \in \{0,1,2,3,4,5,6\}$ and
	$$\widetilde{\sigma}_i (f(A)) = \widetilde{\sigma}_i(f(A_\Lambda)), \qquad f \in \mathcal{M}(A),\quad i \in \{1,2,3,4,5,6\}.
	$$
	Also, $\codim \, \ran(P_\Lambda) < \infty$.
\end{lemma}
\begin{proof}
	The inclusions $\widetilde{\sigma}_i(f(A_\Lambda)) \subseteq \widetilde{\sigma}_i(f(A))$ are given in Lemma \ref{resolventLemma}. Let us show that the inclusions $\sigma_i(f(A)) \subseteq \sigma_i(f(B))$ also hold. To do this, we prove the following claims for all $f \in \mathcal{M}(A)$,
	\begin{enumerate}
		\item If $\nul (f(A_\Lambda)) < \infty$, then $\nul (f(A)) < \infty$.
		\item If $\defec (f(A_\Lambda)) < \infty$, then $\defec (f(A)) < \infty$.
		\item If $\ran(f(A_\Lambda))$ is closed/complemented in $\ran(P_\Lambda)$, then $\ran(f(A))$ is closed/complemented in $X$.
		\item If $\kernel(f(A_\Lambda))$ is complemented in $\ran(P_\Lambda)$, then $\kernel(f(A))$ is complemented in $X$.
	\end{enumerate} 
	Set $\Omega = \widetilde \sigma(A)\setminus \Lambda$, so $\Omega = M_A \setminus \widetilde \sigma_i(A)$ for any $i \in \{1,\ldots,6\}$. Lemma \ref{commutingProjectionsLemma} implies that $\widetilde{\sigma}_i(A|_\Omega) = \emptyset$ for $i \in \{1,2,3,4,5,6\}$, whence $\ran(P_\Omega)$ is finite dimensional, i.e., $\codim\, \ran(P_\Lambda) < \infty$ since $\ran(P_\Lambda)$ and $\ran(P_\Omega)$ are complementary subspaces. Since $\kernel (f(A_\Lambda)) = \kernel(f(A)) \cap \ran(P_\Lambda)$ and $\ran(f(A_\Lambda)) = \ran(f(A)) \cap \ran(P_\Lambda)$ (see the proof of Lemma \ref{resolventLemma}), we conclude that claims (1) and (2) hold true.
	
	For the claim regarding closedness in (3), assume $\ran(f(A_\Lambda))$ is closed in $\ran(P_\Lambda)$, so $\ran(f(A_\Lambda))$ is closed in $X$ too. Since $\ran(f(A)) = \ran(f(A_\Lambda)) \oplus \ran(f(A_\Omega))$, we have that $\ran(f(A))/\ran(f(A_\Lambda))$ is finite-dimensional in $X/\ran(f(A_\Lambda))$, hence closed. Thus $\ran(f(A))$ is closed in $X$.	For the claim regarding complementation in (3), assume $\ran(f(A_\Lambda)) \oplus U = \ran(P_\Lambda)$ for some closed subspace $U$. Note that $\ran(f(A_\Omega)) \oplus V = Z$ for some closed subspace since $\dim \, \ran(P_\Omega) < \infty$. Therefore $\ran(f(A)) \oplus (U \oplus V) = X$, and (3) follows. An analogous reasoning proves the claim (4).
	
	Now, a similar reasoning as above with subspaces $\dom(f(A)), \dom(f(A_\Lambda))$ shows that, if $\infty \in \widetilde{\sigma}_i(f(A))$, then $\infty \in \widetilde{\sigma}_i(f(A_\Lambda))$ for $i \in \{1,2,3,4,5,6\}$. Therefore, $\widetilde{\sigma}_i(f(A)) \subseteq \widetilde{\sigma}_i(f(A_\Lambda))$, as we wanted to show.
	
	Finally, to prove that $M_{A_\Lambda} \subseteq \widetilde{\sigma}_i(A_\Lambda)$ (for $i \in \{1,2,3,4,5,6\}$), note that $M_A \setminus  \widetilde{\sigma}_i(A) = \widetilde{\sigma}(A)\setminus \Lambda \subseteq \rho(A_\Lambda)$ by Lemma \ref{commutingProjectionsLemma}.
\end{proof}

We are now ready to prove the spectral mapping theorems for most of the extended essential spectra considered in the Introduction. For the sake of clarity, we separate the proof of each inclusion into two different propositions.

\begin{proposition}\label{essentialInclusionProp1}
	Let $A \in \BSect(\omega,a)$ and let $f \in \mathcal{M}(A)$ be quasi-regular at $M_A$. Then
	$$\widetilde{\sigma}_i(f(A)) \subseteq f(\widetilde{\sigma}_i(A)), \qquad i \in \{0,1,2,3,4,5,7,8\}.
	$$
\end{proposition}
\begin{proof}
	The inclusion for $\widetilde \sigma_0$ is already given in Proposition \ref{spectralInclusionProp}. For $i \in \{1,2,3,4,5\}$, we can assume $M_A \subseteq \widetilde{\sigma}_i(A)$ without loss of generality by Lemma \ref{quotientNFCLemma}. Thus $\widetilde{\sigma}_i(f(A)) \subseteq f(\widetilde{\sigma}_i(A))$ for $i \in \{1,2,3,4,5\}$ by Proposition \ref{essentialMappingNoSingular}(b).
	
	Now, we show the inclusions for $\widetilde \sigma_7, \widetilde \sigma_8$. So take $i \in \{7,8\}$, and let $\mu \in \CC \setminus f(\widetilde{\sigma}_i(A))$. Note that we can assume $\mu = 0$. If $0 \notin f(M_A)$, Proposition \ref{essentialMappingNoSingular}(b) implies $0 \notin \widetilde \sigma_i(f(A))$. So assume $0 \in f(M_A)$. As $0 \notin f(\widetilde{\sigma}_i(A))$, Lemma \ref{singularPointsEquivalenceLemma} and Remark \ref{isolatedEssentialPoints}(c) imply that $f^{-1}(0)\cap \widetilde{\sigma}(A)$ is a finite set. Let $\lambda_1, ..., \lambda_N$ be an enumeration of $(f^{-1}(0) \cap \widetilde{\sigma}(a) ) \setminus M_A$, and let $n_1,...,n_N$ be the multiplicity of $f$ at $\lambda_1,\ldots,\lambda_N$ respectively. Set $\displaystyle{r(z):= \prod_{j=1}^{N}\left(\frac{\lambda_j - z}{b-z}\right)^{n_j}}$ for some $b \in \rho(A)$, and $g(z): = f(z)/r(z)$ for $z \in \DD$. Several applications of Lemma \ref{fEquivalentgLemma} yield that $g \in \mathcal{M}(A)$ and $f(A) = r(A)g(A) = g(A)r(A)$ with $\dom(f(A)) = \dom(g(A))$. Also, one has $g^{-1}(0) \cap \widetilde{\sigma}(A) \subseteq M_A \setminus \widetilde{\sigma}_i(A)$. Then, $g^{-1}(0) \cap \Lambda = \emptyset$ where $\Lambda := \widetilde{\sigma}(A) \setminus (M_A \setminus \widetilde{\sigma}_i(A))$. Thus, $g(A_\Lambda)$ is invertible by Proposition \ref{spectralInclusionProp} and Lemma \ref{resolventLemma}.
	On the other hand, $\dim\,\ran(I- P_\Lambda) < \infty$ by Lemma \ref{quotientNFCLemma} (see also Lemma \ref{singularPointsEquivalenceLemma}), so $g(A_{\widetilde\sigma(A)\setminus \Lambda}) \in \Phi_i$. As a consequence, $g(A) \in \Phi_i$. Furthermore, $r(A)$ is a bounded operator and belongs to $\Phi_i$. Therefore, $f(A) = r(A)g(A) = g(A) r(A)\in \Phi_i$ by \cite[Theorem I.3.16]{edmunds1987spectral} for $i=7$, and by Lemma \ref{i89Lemma} for $i =8$, that is $0 \notin \widetilde{\sigma}_i(f(A))$.
	
	Now, assume $\infty \notin f(\widetilde{\sigma}_i(A))$. Reasoning as at the end of the proof of Proposition \ref{essentialMappingNoSingular}, we can assume $\rho(f(A)) \neq \emptyset$ without loss of generality. So let $\nu \in \rho(f(A))$. An application of what we have already proven to the function $\displaystyle{\frac{1}{\nu- f(z)}}$ shows that $\displaystyle{0 \notin \widetilde{\sigma_i}\left((\nu - f(A))^{-1}\right)}$, that is $\infty \notin \widetilde{\sigma}_i(f(A))$ for $i \in \{8\}$, and the proof is finished.
\end{proof}

\begin{proposition}\label{essentialInclusionProp2}
	Let $A \in \BSect(\omega,a)$ and let $f \in \mathcal{M}(A)$ be quasi-regular at $M_A$. Then
	$$f(\widetilde{\sigma}_i(A)) \subseteq \widetilde{\sigma}_i(f(A)) , \qquad i \in \{0,1,2,3,4,5,6,8\}.
	$$
\end{proposition}
\begin{proof}
	Note that $\widetilde{\sigma}_6(f(A)) \subset \widetilde{\sigma}_i(f(A))$ for each $i \in \{0,1,2,3,4,5,6,8\}$. Thus, Lemma \ref{singularPointsEquivalenceLemma} yields that it is enough to prove the claim for $i \in \{0,6\}$. Hence, we assume $i \in \{0,6\}$ from now on.
	
	Let $\mu \in f(\widetilde{\sigma}_i(A))$ with $\mu \neq \infty$, so we can assume $\mu = 0$ without loss of generality. If $0 \in  f(\widetilde{\sigma}_i(A)) \setminus f(M_A)$, then $0 \in \widetilde \sigma_i(f(A))$ by Proposition \ref{essentialMappingNoSingular}(a). So assume $0 \in f(\widetilde{\sigma}_i(A))$ with  $0 \in f(M_A)$.  If any point in $f^{-1}(0) \cap \widetilde{\sigma}_i(A)$ is an an accumulation point of $\widetilde{\sigma}_i(A)$ (and we rule out the trivial case where $f$ is constant), then $0$ is an accumulation point of $f(\widetilde{\sigma}_i(A))\setminus f(M_A) \subseteq \widetilde{\sigma}_i(f(A))$ (see Proposition \ref{essentialMappingNoSingular}(a)), thus $0 \in \widetilde{\sigma}_i(f(A))$ since $\sigma_i(T)$ is closed for any $T \in C(X)$. 
	So assume that each point in $f^{-1}(0) \cap \widetilde{\sigma}_i(A)$ is an isolated point in $\widetilde{\sigma}_i(A)$, and set 
	$$V_A := \{d \in f^{-1}(0) \cap \widetilde{\sigma}_i(A) \, | \, d \mbox{ is not an isolated point of } \widetilde{\sigma}(A)\},$$
	which is a finite set by Remark \ref{isolatedEssentialPoints}(c).
	
	Assume first that $V_A$ is not empty (thus $i = 6$). One has that, for each $d \in V_A$, there is some neighbourhood $\Omega_d$ of $d$ such that $\Omega_d \cap \widetilde{\sigma}(A) = \{d, \lambda_1^d, \lambda_2^d,...\}$, where $\lambda_j^d \in \sigma_p(A)\setminus \sigma_i(A)$, each $\lambda_j^d$ is an isolated point of $\sigma(A)$, and $\lambda_j^d \xrightarrow[j\to \infty]{}d$. Thus, $(f^{-1}(0)\cap \widetilde{\sigma}(A))\setminus (\cup_{d \in V_A} \Omega_d)$ is a finite set. Let $\kappa_1, \ldots, \kappa_N$ be the elements of this set, let $n_1,\ldots, n_N$ be the multiplicity of the zero of $f$ at $\kappa_1, \ldots, \kappa_N$ respectively, and set $\displaystyle{g(z) := f(z) \prod_{j=1}^N \left(\frac{b-z}{\kappa_j-z}\right)^{n_j}}$. Several applications of Lemma \ref{domainLemma} yield that $g\in \mathcal{M}(A)$ with $\dom(g(A)) = \dom(f(A))$, and
	\begin{equation}\label{laterforRemark}
		f(A) = \left( \prod_{j=1}^N \left((\kappa_j-A)(b-A)^{-1}\right)^{n_j}\right) g(A) = g(A) \prod_{j=1}^N \left((\kappa_j-A)(b-A)^{-1}\right)^{n_j} ,
	\end{equation}
	where in the last term we regard $(\kappa_j-A)(b-A)^{-1}$ as bounded operators on $\dom(f(A))$. Let us show that $0 \in \widetilde{\sigma}_6(g(A))$, from which follows	 $0 \in \widetilde{\sigma}_6(f(A))$, see for example \cite[Theorem I.3.20]{edmunds1987spectral}. Note that $g^{-1}(0)\cap \widetilde{\sigma}(A) \subset \cup_{d\in V_A} \Omega_d$, which is a countable set. As a consequence, $0$ is an accumulation point of $\CC \setminus g(\widetilde \sigma(A))$. Thus, Proposition \ref{spectralInclusionProp} implies that $0$ is an accumulation point of $\rho(g(A))$. If $0$ is also an accumulation point of $\widetilde{\sigma}(g(A))$, then $0 \in \widetilde{\sigma}_6(g(A))$ by Remark \ref{isolatedEssentialPoints}. So assume that $0$ is not an accumulation point of $\widetilde{\sigma}(g(A))$. Since $\sigma_p(g(A)) \subset g(\sigma_p(A))$ (Lemma \ref{pointSpectraLemma}), and $\lambda_j^d \in \sigma_p(A)$ with $\lambda_j^d \xrightarrow[j \to \infty]{}d$ for each $d \in V_A$, it follows that $g(\lambda_j^d)=0$ for all but finitely many pairs $(j,d) \in \NN \times V_A$. Hence, the set $g^{-1}(0) \cap \sigma_p(A)$ has infinite cardinal, so $\nul(g(A)) \geq \sum_{\lambda \in g^{-1}(0) \cap \sigma_p(A)} \nul(\lambda - A) = \infty$. 
	Then Remark \ref{isolatedEssentialPoints}(c) yields that $0 \in \widetilde{\sigma}_6(g(A))$, as we wanted to prove.

	Now, assume $V_A = \emptyset$, so each $d \in f^{-1}(0) \cap \widetilde{\sigma}_i(A)$ is an isolated point of $\widetilde{\sigma}(A)$. Set $\Lambda := f^{-1}(0) \cap M_A \cap \widetilde \sigma_i(A)$. Note that, in the case $i=6$, then $\dim \ran (P_\Lambda) = \infty$ as a consequence of Lemma \ref{commutingProjectionsLemma}. Since $f(\Lambda) = \{0\}$, we have $\widetilde \sigma (f(A_\Lambda)) \subseteq \{0\}$ by Proposition \ref{spectralInclusionProp}. Hence, $\widetilde\sigma_i(f(A_\Lambda)) = \{0\}$ since $\widetilde{\sigma}_i(f(A_\Lambda))$ cannot be the emptyset (at least for any operator with non-empty resolvent set, see for example \cite{gramsch1971spectral}). Therefore, $0 \in \widetilde{\sigma}_i(f(A))$ by Lemma \ref{resolventLemma}, as we wanted to show.


	Finally, we deal with the case $\mu = \infty$. Reasoning as at the end of the proof of Proposition \ref{essentialMappingNoSingular}, we can assume that $\rho(f(A)) \neq \emptyset$. Take any $\nu \in\rho(f(A))$, so $\infty \in \widetilde{\sigma}_i (f(A))$ if and only if $\displaystyle{0 \in \widetilde{\sigma}_i \left((\nu - f(A))^{-1}\right)}$. But $\displaystyle{0 \in \widetilde{\sigma}_i\left((\nu - f(A))^{-1}\right)}$ by applying what we have already proven to the function $\displaystyle{\frac{1}{\nu - f(z)}}$, and the proof is finished.

\end{proof}


As a consequence, we have the following
\begin{theorem}\label{finalEssentialTheorem}
	Let $A \in \BSect(\omega,a)$ and $f \in \mathcal{M}(A)$ quasi-regular at $M_A$. Then
	\begin{align*}
		\widetilde{\sigma}_i(f(A)) &= f(\widetilde{\sigma}_i(A)), \qquad i \in \{0,1,2,3,4,5,8\},
		\\ f(\widetilde{\sigma}_6(A)) &\subseteq \widetilde{\sigma}_6(f(A)),
		\\ \widetilde{\sigma}_7(f(A)) &\subseteq f(\widetilde{\sigma}_7(A)).
	\end{align*}
\end{theorem}
\begin{proof}
	Immediate consequence of Proposition \ref{essentialInclusionProp1} and Proposition \ref{essentialInclusionProp2}.
\end{proof}

It is known that the spectral mapping theorem does not hold (in general) for $\widetilde{\sigma}_6, \widetilde{\sigma}_7$, see for instance \cite[Section 3]{edmunds1987spectral}. However, we do not know if it holds for $\widetilde{\sigma}_9$ for the regularized functional calculus considered here. Indeed, it holds if $M_A = \emptyset$, see Proposition \ref{essentialMappingNoSingular}.

\section{Final remarks}\label{finalRemarksSection}

\subsection{Operators with bounded functional calculus}
A natural question is whether the condition of quasi-regularity can be relaxed in Theorem \ref{finalEssentialTheorem}. A possible candidate for this relaxed condition could be asking for $f$ to have well-defined limits at $M_A$. In order to prove a spectral mapping theorem for this wider class of functions we will need to impose the following condition on $A$. The property of having a bounded (regularized) functional calculus is studied in , we ask the operator $A$ to fulfil the following property, which is studied in \cite{cowling1996banach, haase2006functional, morris2010local}.

\begin{definition}
	Let $A \in \BSect(\omega,a)$. We say that the regularized functional calculus of $A$ is bounded if $f(A) \in \mathcal{L}(X)$ for every bounded $f \in \mathcal{M}(A)$.
\end{definition}

\begin{lemma}\label{notInjectivityLemma}
	Let $A \in \BSect(\omega,a)$ and $f \in \mathcal{M}(A)$. Then $f$ is regular at $\sigma_p(A) \cap M_A$.
\end{lemma}
\begin{proof}
	The proof is analogous to the case of sectorial operators, see \cite[Lemma 4.2]{haase2005general}.
\end{proof}

\begin{theorem}\label{boundedFunctionalCTheorem}
	Let $A \in \BSect(\omega,a)$ such that the regularized functional calculus of $A$ is bounded, and let $f \in \mathcal{M}(A)$ with (possibly $\infty$-valued) limits at $M_A$. Then
	\begin{align*}
		\widetilde{\sigma}_i(f(A)) &= f(\widetilde{\sigma}_i(A)), \qquad i \in \{0,1,2,3,4,5,8\},
		\\ f(\widetilde{\sigma}_6(A)) &\subseteq \widetilde{\sigma}_6(f(A)),
		\\ \widetilde{\sigma}_7(f(A)) &\subseteq f(\widetilde{\sigma}_7(A)).
	\end{align*}
\end{theorem}
\begin{proof}
	The proof of this claim is completely analogous to the path followed in this paper to prove Theorem \ref{finalEssentialTheorem}. Indeed, the quasi-regularity notion is only explicitly needed in the proofs of Lemma \ref{differentNFCLemma} and Proposition \ref{spectralInclusionProp}. This creates a `cascade' effect, and all following results need the quasi-regularity assumption in order to apply Proposition \ref{spectralInclusionProp}. Therefore, we will prove the claim once we prove the following version of Proposition \ref{spectralInclusionProp}:
	
	``Let $A \in \BSect(\omega,a)$ such that the regularized functional calculus of $A$ is bounded, and let $f \in \mathcal{M}(A)$ with (possibly $\infty$-valued) limits at $M_A$. Then $\widetilde{\sigma}(f(A)) \subset f(\widetilde{\sigma}(A))$.''
	
	We outline  the proof of this claim. Let $\mu \in \CC_\infty$ with $\mu \notin f(\widetilde{\sigma})$, and set $\displaystyle{f_\mu} = \frac{1}{\mu - f}$ if $\mu \in \CC$ or $f_\mu = f$ if $\mu = \infty$, and we will show that $f_\mu \in \mathcal{M}(A)$ with $f_\mu(A)  \in \mathcal{L}(X)$. Note that $f_\mu$ has finite limits at $M_A$. Even more, $f_\mu$ is regular at $\sigma_p(A) \cap M_A$ by Lemma \ref{notInjectivityLemma}. Proceeding as in the proof of Lemma \ref{differentNFCLemma}, we can assume that $f_\mu$ has finitely many poles, all of them contained in $\rho(A)$. Let  $\displaystyle{r(z) := \prod_{j=1}^{n}\frac{(\lambda_j-z)^{n_j}}{(b-z)^{n_j}}}$, where $\lambda_j, n_j$ are the poles of $f_\mu$ and their order, respectively. Hence, $r f_\mu$ has no poles, is regular at $\sigma_p(A) \cap M_A$ and has finite limits at $M_A$, thus $r f_\mu$ is bounded. For any $b \in \rho(A)$, the function $\displaystyle{h(z):=\frac{1}{b-z} \prod_{d\in \{-a,a\}\setminus \sigma_p(A)} \frac{z-a}{b-z}}$ regularizes $rf_\mu$, so $rf_\mu \in \mathcal{M}_A$. Since the regularized functional calculus of $A$ is bounded, then $rf_\mu (A) \in \mathcal{L}(X)$. Moreover, $r(A)$ is bounded and invertible. Therefore, $hr$ regularizes $f_\mu$ with $f_\mu(A) = r(A)^{-1} (rf_\mu)(A) \in \mathcal{L}(X)$, and the claim follows.
\end{proof}

\subsection{Spectral mapping theorems for sectorial operators and strip-type operators}\label{sectorialSubSect}

Lemmas \ref{functionalCalculusProperties}, \ref{pointSpectraLemma}, \ref{regularizerLemma} and \ref{differentNFCLemma} are the only properties of the regularized functional calculus of bisectorial-like operators that are used in the proofs given in Sections \ref{noSingularPointsSection} and \ref{spectralMappingSection}. Therefore, for any other regularized functional calculus satisfying such properties, one can prove spectral mapping theorems for essential spectra analogous to the ones given in Theorem \ref{finalEssentialTheorem}.
In particular, one has the following result for the regularized functional calculus of sectorial operators and the regularized functional calculus of strip operators considered in \cite{haase2005spectral}.

\begin{theorem}\label{spectralMappingThSectorial}
	Let $A$ be a sectorial operator of angle $\phi \in [0,2\pi)$, and let $f$ be a function in the domain of the regularized functional calculus of $A$ that is quasi-regular at $\{0,\infty\}\cap \widetilde{\sigma}(A)$. Then
	\begin{align*}
		\widetilde{\sigma}_i(f(A)) &= f(\widetilde{\sigma}_i(A)), \qquad i \in \{0,1,2,3,4,5,8\},
		\\ f(\widetilde{\sigma}_6(A)) &\subseteq \widetilde{\sigma}_6(f(A)),
		\\ \widetilde{\sigma}_7(f(A)) &\subseteq f(\widetilde{\sigma}_7(A)).
	\end{align*}
\end{theorem}

\begin{theorem}\label{spectralMappingThStrip}
	Let $A$ be a strip operator of height $h \geq 0$, and let $f$ be a function in the domain of the regularized functional calculus of $A$ that is quasi-regular at $\{\infty\}\cap \widetilde{\sigma}(A)$. Then
	\begin{align*}
		\widetilde{\sigma}_i(f(A)) &= f(\widetilde{\sigma}_i(A)), \qquad i \in \{0,1,2,3,4,5,8\},
		\\ f(\widetilde{\sigma}_6(A)) &\subseteq \widetilde{\sigma}_6(f(A)),
		\\ \widetilde{\sigma}_7(f(A)) &\subseteq f(\widetilde{\sigma}_7(A)).
	\end{align*}
\end{theorem}
\subsection{A spectral mapping theorem for the point spectrum}\label{pointSubsection}
To finish this paper, we give a spectral mapping theorem for the point spectrum. To prove it, we need to restrict to functions $f \in \mathcal M_A$ satisfying the following condition: \medskip

\noindent (\textbf{P}) \quad For each $d \in M_A$ such that $f(d) \notin f(\sigma_p(A)) \cup \{\infty\}$, there is some $\beta>0$ for which
\begin{itemize}
	\item if $d \in \CC$, then $|f(z) - c_d| \gtrsim |z-d|^\beta$ as $z\to d$, or
	\item if $d=\infty$, then $|f(z) - c_d| \gtrsim |z|^{-\beta}$ as $z\to d$,
\end{itemize}
where $c_d$ denotes the limit of $f(z)$ as $z \to d$. \medskip

\begin{proposition}\label{PointspectralMappingPropositionAddition}
	Let $A \in \BSect(\omega,a)$ and $f \in \mathcal{M}(A)$ such that $f$ is quasi-regular at $M_A$. Then
	\begin{align*}
		f(\sigma_p(A)) & \subseteq \sigma_p (f(A)) \subseteq f(\sigma_p(A)) \cup f(M_A).
	\end{align*}	
	If, furthermore, $f$ satisfies condition {\normalfont(\textbf{P})} above, then $f(\sigma_p(A)) = \sigma_p (f(A))$.
\end{proposition}
\begin{proof}
	The proof of the inclusions $f(\sigma_p(A))  \subset \sigma_p (f(A)) \subset f(\sigma_p(A)) \cup f(M_A)$ runs the same as for sectorial operators, see \cite[Proposition 6.5]{haase2005spectral}. Assume now that $f$ satisfies (\textbf{P}). All that is left to prove is that, if $\mu \in f(M_A) \backslash f(\sigma_p(A))$, then $\mu \notin \sigma_p(f(A))$. The statement is trivial if $\mu = \infty$, so assume $\mu \in \CC \backslash f(\sigma_p(A))$, and consider the function $g:=\frac{1}{\mu-f}$, which is quasi-regular at $M_A$. Note that poles of $g$ are precisely $f^{-1}(\mu) \subset \CC \backslash \sigma_p(A)$. Moreover, $g$ is regular at $M_A \cap \sigma_p(A)$, since by assumption $\mu \notin f(\sigma_p(A))$. Let now
	\begin{align*}
		h_{l,m,n}(z):= \frac{(z-a)^m (z+a)^n}{(b-z)^{l+m+n}}, \qquad z \in \CC, \quad l,m,n \in \NN, \quad \, b>a.
	\end{align*}
	Then, by the assumptions made on $f$, $h_{l,m,n}g$ is regular for some $m,n,l$ large enough, and where $l, m, n =0$ if $\infty, a, -a \notin \widetilde{\sigma}(A) \backslash \sigma_p(A)$ respectively. Since $h_{l,m,n}(A) = (A-a)^m(A+a)^n R(b,A)^{m+n+l}$ is bounded and injective, $h_{l,m,n}$ regularizes $g$. Hence $g \in \mathcal{M}(A)$, which by Lemma \ref{functionalCalculusProperties}(4) implies that $\mu-f$ is injective, as we wanted to show.
\end{proof}

\section*{Acknowledgements}
This research has been partially supported by Project PID2019-105979GB-I00 of the MICINN, and by BES-2017-081552, MINECO, Spain and FSE.

\bibliography{manuscript}

\providecommand{\bysame}{\leavevmode\hbox to3em{\hrulefill}\thinspace}
\providecommand{\MR}{\relax\ifhmode\unskip\space\fi MR }
\providecommand{\MRhref}[2]{%
  \href{http://www.ams.org/mathscinet-getitem?mr=#1}{#2}
}
\providecommand{\href}[2]{#2}
\begin{thebibliography}{10}

\bibitem{arendt2006maximal}
W.~Arendt and M.~Duelli, \emph{Maximal $l^p$-regularity for parabolic and
  elliptic equations on the line}, J. Evol. Equ. \textbf{6} (2006), no.~4, 773.

\bibitem{bade1953operational}
W.~Bade, \emph{An operational calculus for operators with spectrum in a strip},
  Pac. J. Math. \textbf{3} (1953), no.~2, 257--290.

\bibitem{browder1961spectral}
F.E. Browder, \emph{On the spectral theory of elliptic differential operators.
  {I}}, Math. Ann. \textbf{142} (1961), no.~1, 22--130.

\bibitem{cowling1996banach}
M.~Cowling, I.~Doust, A.~Mcintosh, and A.~Yagi, \emph{Banach space operators
  with a bounded {$H^\infty$} functional calculus}, J. Aust. Math. Soc.
  \textbf{60} (1996), no.~1, 51--89.

\bibitem{dunford1963linear}
N.~Dunford and J.T. Schwartz, \emph{Linear operators. {II}}, Interscience
  Publishers, New York, 1963.

\bibitem{edmunds1987spectral}
D.E. Edmunds and W.D. Evans, \emph{Spectral theory and differential operators},
  Oxford University Press, Oxford/{New York}, 1987.

\bibitem{gonzalez1985spectral}
M.~Gonz{\'a}lez and V.M. Onieva, \emph{On the spectral mapping theorem for
  essential spectra}, Publ. Sec. Mat. Univ. Aut{\`o}noma Barcelona (1985),
  105--110.

\bibitem{gramsch1971spectral}
B.~Gramsch and D.~Lay, \emph{Spectral mapping theorems for essential spectra},
  Math. Ann. \textbf{192} (1971), no.~1, 17--32.

\bibitem{gustafson1969essential}
K.~Gustafson and J.~Weidmann, \emph{On the essential spectrum}, J. Math. Anal.
  Appl. \textbf{25} (1969), no.~1, 121--127.

\bibitem{haase2005general}
M.~Haase, \emph{A general framework for holomorphic functional calculi}, Proc.
  Edinb. Math. Soc. \textbf{48} (2005), no.~2, 423--444.

\bibitem{haase2005spectral}
\bysame, \emph{Spectral mapping theorems for holomorphic functional calculi},
  J. London Math. Soc. \textbf{71} (2005), no.~3, 723--739.

\bibitem{haase2006functional}
\bysame, \emph{The functional calculus for sectorial operators}, vol. 169,
  Oper. Theory Adv. Appl., Birkh{\"a}user, Basel, 2006.

\bibitem{kato1966perturbation}
T.~Kato, \emph{Perturbation theory for linear operators}, Springer-Verlag,
  Berlin, 1966.

\bibitem{lay1968characterizations}
D.~Lay, \emph{Characterizations of the essential spectrum of {FE Browder}},
  Bull. Amer. Math. Soc. \textbf{74} (1968), no.~2, 246--248.

\bibitem{mcintosh1986operators}
A.~McIntosh, \emph{Operators which have an $h_\infty$ functional calculus},
  Miniconference on operator theory and partial differential equations,
  vol.~14, Austral. Nat. Univ., Mathematical Sciences Institute, Canberra,
  1986, pp.~210--232.

\bibitem{mielke1987maximalel}
A.~Mielke, \emph{{\"U}ber maximalel p-regularit{\"a}t f{\"u}r
  differentialgleichungen in banach-und hilbert-r{\"a}umen}, Math. Ann.
  \textbf{277} (1987), no.~1, 121--133.

\bibitem{morris2010local}
A.J. Morris, \emph{Local quadratic estimates and the holomorphic functional
  calculi}, The AMSI--ANU Workshop on Spectral Theory and Harmonic Analysis,
  vol.~44, Proc. Centre Math. Appl. Austral. Nat. Univ., 2010, pp.~211--231.

\bibitem{pietsch1960theorie}
A.~Pietsch, \emph{Zur {Theorie} der {$\sigma$-Transformationen} in
  lokalkonvexen {Vektorr{\"a}umen}}, Math. Nach. \textbf{21} (1960), no.~6,
  347--369.

\bibitem{schechter1966essential}
M.~Schechter, \emph{On the essential spectrum of an arbitrary operator. {I}},
  J. Math. Anal. Appl. \textbf{13} (1966), no.~2, 205--215.

\bibitem{taylor1958introduction}
A.E. Taylor and D.C. Lay, \emph{Introduction to {Functional Analysis}}, Wiley
  \& sons, New York, 1958.

\bibitem{wolf1959essential}
F.~Wolf, \emph{On the essential spectrum of partial differential boundary
  problems}, Commun. Pure Appl. Math. \textbf{12} (1959), no.~2, 211--228.

\bibitem{yood1951properties}
B.~Yood, \emph{Properties of linear transformations preserved under addition of
  a completely continuous transformation}, Duke Math. J. \textbf{18} (1951),
  no.~3, 599--612.

\end{thebibliography}
\bibliographystyle{amsplain}

\end{document}